\documentclass[10pt]{article}
\usepackage{amssymb}
\usepackage{graphicx}
\usepackage{xcolor} 
\usepackage{tensor}
\usepackage{fullpage} 
\usepackage{amsmath}
\usepackage{amsthm}
\usepackage{verbatim}

\usepackage{enumitem}
\setlist[enumerate]{leftmargin=1.5em}
\setlist[itemize]{leftmargin=1.5em}

\usepackage{seqsplit}

\definecolor{green}{rgb}{0,0.8,0} 

\newtheorem{maintheorem}{Theorem}

\newtheorem{theorem}{Theorem}[section]

\newtheorem{lemma}[theorem]{Lemma}
\newtheorem{proposition}[theorem]{Proposition}

\theoremstyle{definition}
\newtheorem{definition}[theorem]{Definition}

\theoremstyle{remark}
\newtheorem{remark}[theorem]{Remark}

\numberwithin{equation}{section}
\newcommand{\nrm}[1]{\Vert#1\Vert}

\newcommand{\nnrm}[1]{{\vert\kern-0.25ex\vert\kern-0.25ex\vert #1 
    \vert\kern-0.25ex\vert\kern-0.25ex\vert}}

\newcommand{\lap}{\Delta}

\newcommand{\rd}{\partial}
\newcommand{\nb}{\nabla}

\newcommand{\alp}{\alpha}

\newcommand{\gmm}{\gamma}

\newcommand{\eps}{\epsilon}

\newcommand{\Lmb}{\Lambda}

\newcommand{\tht}{\theta}

\newcommand{\Omg}{\Omega}

\newcommand{\bbR}{\mathbb R}

\newcommand{\bbT}{\mathbb T}

\newcommand{\tu}{\tilde{u}}					

\newcommand{\trho}{\tilde{\rho}}

\newcommand{\R}{\mathbb R}

\newcommand{\pa}{\partial}

\newcommand{\e}{\varepsilon}
\newcommand{\lt}{\left}
\newcommand{\rt}{\right}
\newcommand{\bq}{\begin{equation}}
\newcommand{\eq}{\end{equation}}

\begin{document}

\title{On well-posedness and singularity formation \\
	 for the Euler--Riesz system}
\author{Young-Pil Choi\thanks{Department of Mathematics, Yonsei University, Seoul 03722, Republic of Korea. E-mail: ypchoi@yonsei.ac.kr} \and In-Jee Jeong\thanks{Department of Mathematics, Korea Institute for Advanced Study, Seoul 02455, Republic of Korea. E-mail: ijeong@kias.re.kr}}
\date{\today}



\maketitle


\begin{abstract}
In this paper, we investigate the initial value problem for the Euler--Riesz system, where the interaction forcing is given by $\nb(-\lap)^{s}\rho$ for some $-1<s<0$, with $s = -1$ corresponding to the classical Euler--Poisson system. We develop a  functional framework to establish local-in-time existence and uniqueness of classical solutions for the Euler--Riesz system. In this framework, the fluid density could decay fast at infinity, and the Euler--Poisson system can be covered as a special case. Moreover, we prove local well-posedness for the pressureless Euler--Riesz system when the potential is repulsive, by observing hyperbolic nature of the system. Finally, we present sufficient conditions on the finite-time blowup of classical solutions for the isentropic/isothermal Euler--Riesz system with either attractive or repulsive interaction forces. The proof, which is based on estimates of several physical quantities, establishes finite-time blowup for a large class of initial data; in particular, it is not required that the density is of compact support. 
\end{abstract}


\section{Introduction}

\subsection{Derivation of the system}

This paper is devoted to the study of the well-posedness theory for the \textit{Euler--Riesz} system
\begin{equation}  \label{eq:Euler-second}
\left\{
\begin{aligned}
&\rd_t \rho + \nb_x\cdot (\rho u) = 0, \\
&\rd_t (\rho u) + \nabla_x\cdot ( \rho u\otimes u) + c_p \nb_x(\rho^\gmm )=  c_K \rho \nb_x\Lmb^{\alp-d}\rho
\end{aligned}
\right.
\end{equation} either in $\Omg = \bbR^d$ or $\bbT^d$, where $\rho(t,\cdot) : \Omg \rightarrow \bbR_{+} $ and $u(t,\cdot) : \Omg \rightarrow \bbR^d$ denote the density and the velocity of the fluid, respectively. Here, the Riesz operator $\Lmb^s$ is defined by $(-\Delta_x)^{\frac{s}{2}}$ and we shall consider the range $d-2 < \alp < d$. We take $c_p\ge 0$ and $\gmm\ge 1$, with $c_p = 0$ corresponding to the pressureless Euler--Riesz system. The case $c_K>0$ and $c_K <0$ correspond to attractive and repulsive potentials, respectively. In the following, we shall always assume that $\rho>0$, which allows to rewrite the equation for $\rho u$ in \eqref{eq:Euler-second} as 
\[
\rd_t u + u\cdot\nb_x u + c_p \gmm\rho^{\gmm-2}\nb_x \rho = c_K\nb_x \Lmb^{\alp-d}\rho. 
\]


The Euler--Riesz system \eqref{eq:Euler-second} can be derived from the second-order particle system corresponding to Newton's law. More precisely, let $x_i(t)$ and $v_i(t)$ be the position and velocity of the $i$-th particle at time $t >0$. Then the dynamics of interacting point particles through the interaction force $\nabla_x K$ can be described by
\begin{align}\label{par_sys}
\left\{
\begin{aligned}
\frac{d x_i(t)}{dt} &= v_i(t),\quad i=1,\dots,N, \quad t > 0,\cr
\frac{d v_i(t)}{dt} &= \frac{c_K}N \sum_{j \neq i} \nabla_x K (x_i(t) - x_j(t)),
\end{aligned}
\right.
\end{align}
where the interaction potential $K$ is given by the Riesz kernel\footnote{$K\star \rho$ can be given as an inverse fractional Laplacian with an appropriate positive constant $c_{\alpha,d}$, i.e., \[
c_{\alpha,d}K\star \rho = (-\Delta_x)^{\frac{\alpha-d}{2}} \rho =\Lmb^{\alpha-d} \rho .
\]
However, for simplicity, we set $c_{\alpha,d}=1$ throughout this paper.}:
\bq\label{ker}
K(x) = |x|^{-\alpha}, \quad d-2 < \alpha < d.
\eq
Note that the case $\alpha = d-2$ corresponds to the Coulombian kernel when $d\geq 3$, which is the classical model for for the plasma dynamics when $c_K <0$ and for astrophysics problems when $c_K >0$.

When the number of particles $N$ is sufficiently large, it is reasonable to consider the corresponding continuum model in order to reduce the complexity of the particle system, which is important in any practical applications. The classical way is to introduce the number density function $f = f(t,x,v)$ in the phase space $(x,v) \in \R^d \times \R^d$. To be more specific, at the formal level, as the number of $N$ goes to infinity, by means of BBGKY hierarchies or mean-field limits \cite{CCHS19,CCS19,Ha14,HJ07,HJ15,HS19,JW18,LP17,Sepre}, we can derive the following Vlasov equation with Riesz interaction forces from the particle system \eqref{par_sys}:
\bq\label{kin_sys}
\pa_t f + v \cdot \nabla_x f + (c_K\nabla_x K \star \rho) \cdot \nabla_v f = 0, \quad (x,v) \in \R^d, \quad t > 0,
\eq
where $\rho = \rho(t,x)$ is the local particle density given by $\rho = \int_{\R^d} f\,dv$. For the Coulomb kernal, which corresponds to the case $\alpha = d-2$ in \eqref{ker} when $d\geq 3$, the kinetic equation \eqref{kin_sys} is called the Vlasov--Poisson equation, and for $\alpha = d-1$ it is called the Vlasov--Manev equation \cite{IVDB98}. We refer to \cite{Gl96,Pe04} for the general kinetic theory. 

Since the kinetic model \eqref{par_sys} is posed in $2d+1$ dimensions, it is still computationally expensive. For that reason, the kinetic equations are usually being reduced to hydrodynamic equations by suitable asymptotic limits. We consider the local moment estimates in $v$: multiplying the equation \eqref{kin_sys} by $1$ and $v$, respectively, and integrating the resulting equations with respect to $v$, we can derive the following local conservation laws\footnote{One may also derive the equation for the local energy $\rho E = \int_{\R^d} |v|^2 f\,dv$, however, it is not necessarily required in this formal derivation.}:
\begin{align}\label{hydro_sys}
\left\{
\begin{aligned}
&\pa_t \rho + \nabla_x \cdot (\rho u) =0, \quad x \in \R^d, \quad t > 0,\cr
&\pa_t (\rho u) + \nabla_x \cdot (\rho u \otimes u) + \nabla_x \cdot \lt( \int_{\R^d} (v - u) \otimes (v-u)f\,dv \rt)=  c_K\rho \nabla_x K \star \rho.
\end{aligned}
\right.
\end{align}
Here, $u$ denotes the local particle velocity given by $u = \int_{\R^d} vf\,dv/\rho$.
Although the momentum equations in the above are not closed due to the third term on the left hand side, we can formally take the \textit{mono-kinetic} $f(x,v) = \rho(x) \otimes \delta_{u(x)}(v)$ or the local Maxwellian ansatz $f(x,v) = \rho(x) \exp (- |v-u(x)|^2/2)$ in \eqref{hydro_sys}. These different choices of ans\"{a}tze lead to the pressureless and the isothermal Euler--Riesz system, respectively. Recently, rigorous derivations of hydrodynamic models from kinetic models have been studied in \cite{CC20,CCJpre,FK19,Ka18,KMT15}. The mono-kinetic ansatz can be justified by adding the strong local alignment force $(1/\e) \nabla_v \cdot ((v-u)f)$ on the right hand side of \eqref{kin_sys} and considering the limit $\e \to 0$. On the other hand, the isothermal Euler-type system can be rigorously derived by taking into account the strong local alignment and diffusive forces, $(1/\e) \nabla_v \cdot (\nabla_v f + (v-u)f)$. More recently, the rigorous derivation of the pressureless Euler--Riesz system has been established in \cite{Sepre}. Under certain regularity assumptions on solutions to \eqref{eq:Euler-second}, the system \eqref{eq:Euler-second} with $c_p = 0$ can be rigorously derived from the particle system \eqref{par_sys} in the case of mono-kinetic data. We also refer to \cite{CCpre} for the rigorous derivation of hydrodynamic collective behavior models from many interacting particle systems. 

Here we introduce several notations used throughout this paper. For simplicity, we omit $x$-dependence of differential operators, for instances, $\pa f:= \pa_x f$, $\nabla f:= \nabla_x f$, and $\Delta f := \Delta_x f$, etc. We denote by $C$ a generic positive constant and it may differ from line to line. Finally, $f \lesssim_{\alpha,\beta,\cdots} g$ represents that there exists a positive constant $C=C(\alpha,\beta,\cdots ) > 0$ such that $f \leq Cg$.

\subsection{Main Results}

The first main result of this work is the local well-posedenss of classical solutions for the Euler--Riesz system. In particular, derivation of \eqref{eq:Euler-second} from the particle dynamics is justified \cite[Appendix A]{Sepre}. The study of the well-posedness theory for the Euler--Poisson system with either attractive or repulsive force, i.e., the case $\alpha = d-2$ is by now well-established. The existence theory and sticky dynamics of solutions for the pressureless Euler-type system are discussed in \cite{BGSW13,BG98,ERS96,NS09,NT08}. For the Euler--Poisson system, the local/global-in-time existence of strong solutions is studied in  \cite{BK18,Gu98,GP11,PRV95}. In particular, in a recent work \cite{BK18}, the local-in-time existence and uniqueness of strong solutions for the Euler--Poisson system with $\gamma \in (1,5/3)$ in three dimensions are obtained in a weighted Sobolev spaces of fractional order. In particular, in that work, the constructed local-in-time solutions can contain vacuum in some region, i.e., there is no strictly positive lower bound for the density. We point out that although the functional framework developed in this work cannot treat vacuum (region where $\rho$ vanishes), it does cover the Euler--Poisson case without any difficulty, and therefore provides an alternative way of proving well-posedness of classical solutions to the Euler--Poisson system with density decaying fast at infinity. 

Let us state our main well-posedness result in a somewhat rough manner; the precise statements are given in Theorems \ref{thm:pressureless-lwp} (pressureless and repulsive case) and \ref{thm:lwp-pressure} (pressure case), respectively. 
\begin{maintheorem}[Local well-posedness] \label{thm:main1} 
	Consider the system \eqref{eq:Euler-second} on $\Omg = \bbR^d$ or $\bbT^d$ with $d\ge 1$ and $d-2\le \alp < d$. We have the following local well-posedness results. 
	\begin{enumerate}
		\item (pressureless and repulsive case) In the case $c_p = 0$ and $c_K < 0$, for any $m> \frac{d}{2}+2$ there exists a space $X^m(\Omg)$ such that the Euler--Riesz system is locally well-posed: for any $(u_0,\rho_0) \in X^m$, there exists $T>0$ and a unique solution $(u,\rho)$ satisfying the initial condition and belonging to $L^\infty([0,T);X^m)$. In the case $\Omg = \bbT^d$, $X^m$ can be identified with a Sobolev space; we can take $\nrm{(u,\rho)}_{X^m} = \nrm{\rho}_{H^m} + \nrm{u}_{H^{m+\frac{d-\alp}{2}}}$. 
		\item (pressure case) When $c_p>0$ and $1 \le \gmm \le \frac{5}{3}$, there is a space $Y^m(\Omg)$ for any $m> \frac{d}{2}+1$ such that the Euler--Riesz system  with any $c_K \in \bbR$ is locally well-posed: for any $(u_0,\rho_0) \in Y^m$, there exists $T>0$ and a unique solution $(u,\rho)$ satisfying the initial condition and belonging to $L^\infty([0,T);Y^m)$. The $Y^m$-norm is given by $\nrm{u}_{H^m} + \nrm{\rho^{\frac{\gmm-1}{2}}}_{H^m} + \nrm{\nb\ln\rho}_{L^\infty}$ for $\gmm>1$ and in the case $\gmm = 1$, $\rho^{\frac{\gmm-1}{2}}$ should be replaced with $\ln\rho$. 
	\end{enumerate} 
\end{maintheorem}

A detailed discussion involving the statements will be given below in Subsection \ref{subsec:ideas} and further in Section \ref{sec:lin}. For now, let us just briefly comment on the spaces used in Theorem \ref{thm:main1}: in the pressureless case, we can prove local well-posedness simply using the norm $\nrm{\rho}_{H^m} + \nrm{u}_{H^{m+\frac{d-\alp}{2}}}$ when either the domain is bounded (e.g. $\bbT^d$) or $\rho$ attains a uniform positive lower bound in space. Possible decay of $\rho$ at spatial infinity brings some technical difficulties, which is handled by suitably modifying the $H^m$ norm for $\rho$, in a way that additional decay of the derivatives for $\rho$ is encoded. In the case with pressure, the quantity $\nrm{\nb\ln\rho}_{L^\infty}$ plays a similar role. 

Our second main result concerns finite-time singularity formation for the Euler--Riesz system, within the framework of local well-posedness established in Theorem \ref{thm:main1}. Finite-time blow-up of smooth solutions for the compressible Euler or Navier--Stokes system is shown in \cite{CH1, LPZ,S85, Xi} based on the finite propagation speed of the supports of solutions. For the one dimensional Euler--Poisson system, the singularity formation is discussed in \cite{En96}. Critical thresholds between the subcritical region with global-in-time regularity and the supercritical region with finite-time blowup of classical solutions are also analyzed in one dimension \cite{CCZ16,ELT01}. For the multi-dimensional isentropic Euler--Poisson system without radial symmetry, a priori estimate of the finite-time blow-up of solutions is obtained in \cite{W14}. We refer to \cite{CT1, Ch05} and references therein for the general survey on the Euler-type equations and related mathematical problems. 

In the theorem below, we briefly state our second main result on the finite-time singularity formation for the system \eqref{eq:Euler-second} in the presence of pressure. For the precise statements, see Theorems \ref{thm_att} and \ref{thm_rep} (isentropic pressure case) and \ref{thm_iso} (isothermal pressure case) in Section \ref{sec:sing}.
\begin{maintheorem}[Singularity formation] \label{thm:main2} We consider \eqref{eq:Euler-second} in the following cases. 
\begin{enumerate}
\item (isentropic pressure case: $\gamma > 1$) In the attractive case ($c_K > 0$), we assume $\alpha \geq \max\{2, d(\gamma-1)\}$. On the other hand, in the repulsive case ($c_K < 0$), we suppose
\[
1 < \gamma \leq 1 + \frac{2}{d} \quad \mbox{and} \quad \alpha \geq d(\gamma-1).
\]
\item (isothermal pressure case: $\gamma=1$) In the attractive case, we assume $\alpha \geq 2$. We do not require any restriction on $\alpha$ in the repulsive case. 
\end{enumerate}
Then, in both cases, there exist initial data $(\rho_0, u_0)$ which generate a finite-time singularity for the system \eqref{eq:Euler-second}.
\end{maintheorem}

In the attractive case, our strategy does not require any restriction on $\gamma \geq 1$. On the other hand, in the repulsive and isentropic pressure case, we need the condition $\gamma \leq 1 + \frac2d$. However, in this case the restriction $\gamma \leq \frac53$ in Theorem \ref{thm:main1}.2 can be dropped (see Section \ref{subsec:ideas} and Remark \ref{rem:pressure}). Thus, restrictions on $\alpha$ and $\gamma$ for singularity formation are stronger than those for local well-posedness. It is also worth noticing that $\alpha>0$ is required for the proof of singularity formation in the isentropic pressure case, but our proof of blow-up in the isothermal case covers the range $\alpha<0$.

\subsection{Ideas of the proof}\label{subsec:ideas}

Let us give some technical ideas involved in the proofs of main results.

\bigskip

\noindent \textbf{Local well-posedness for pressureless case} (Theorem \ref{thm:pressureless-lwp}). Note that when the pressure term is absent ($c_p=0$), local well-posedness does not follow from a simple energy estimate for $\rho$ and $u$. This can be seen as follows: propagating the $H^s$-norm of $\rho$ requires $\nb\cdot u \in H^s$, while propagating $\nb\cdot u$ in $H^s$ requires $\Lmb^{\alp-d}\lap\rho \in H^s$. Since $\alp-d>-2$, we do not have $\nrm{\Lmb^{\alp-d}\lap\rho}_{H^s} \lesssim \nrm{\rho}_{H^s}$, unlike the Euler--Poisson case where $\alp-d = -2$. To overcome this difficulty, one needs essentially rewrite the system in terms of $\Lmb^{\frac{d-\alp}{2}}u$ and $\rho$. This choice is motivated from linear analysis presented in Section \ref{sec:lin}. Indeed, in the nonlinear system, one can show that there is a cancellation of top order terms in the sum $\nrm{ \rho^{\frac{1}{2}}|\nb|^{s+\frac{d-\alp}{2}} u }_{L^2}^2 + \nrm{|\nb|^s\rho}_{L^2}^2$, when the potential is \textit{repulsive} ($c_K<0$), which allows for the estimate \begin{equation}\label{eq:est}
\begin{split}
\left|\frac{d}{dt}\int_{\R^d} \rho | |\nb|^{s+\frac{d-\alp}{2}} u|^2 -c_K ||\nb|^s \rho|^2\,dx \right| \le C ( \nrm{u}_{H^{s+\frac{d-\alp}{2}}}  + \nrm{\rho}_{H^s})^3.  
\end{split}
\end{equation} On the contrary, when the potential is attractive, the system must be ill-posed. This is hinted in the linear analysis provided in Section \ref{sec:lin}, and we expect that a rigorous proof of nonlinear ill-posedness could be given along the lines of \cite{KVW}. Now, in the proof of \eqref{eq:est}, the key tool is a sharp commutator estimate for $[\Lmb^{\eps},v\cdot\nb]$ which was derived in \cite{CCCGW}; for $v$ smooth, we have \begin{equation*}
\begin{split}
\nrm{ [\Lmb^{\eps},v\cdot\nb]f }_{L^2} \lesssim \nrm{f}_{H^{\eps}},
\end{split}
\end{equation*} see Lemma \ref{lem:comm-Lmb} for a precise statement. The authors in \cite{CCCGW} used this estimate to prove local well-posedness of the active scalar system \begin{equation}\label{eq:sqg-singular}
\begin{split}
\rd_t\tht + \nb^\perp \Lmb^{-\alp}[\tht] \cdot \nb \tht = 0
\end{split}
\end{equation} defined for $\tht: [0,T]\times\bbT^2\rightarrow \bbR$ for some $\alp>0$. Indeed, the equation \eqref{eq:sqg-singular} has some similarities with our pressureless and repulsive system. In this regard, we emphasize that the \textit{parity} of the symbol is crucial in well-posedness for these models; when $\nb^\perp\Lmb^{-\alp}$ in \eqref{eq:sqg-singular} is replaced with some even symbol of the same order, the equation can be shown to be \textit{ill-posed} in Sobolev spaces \cite{KVW}. In the case of \eqref{eq:Euler-second}, when $c_p=0$ and $\nb\Lmb^{\alp-d}\rho$ is replaced with $\Lmb^{1+\alp-d}$ (which is of the same order), the system becomes ill-posed; see linear analysis in the following section. Returning to \eqref{eq:est}, we can close an a priori estimate immediately if $\rho>0$ attains a uniform positive lower bound, since then $\nrm{ \rho^{\frac{1}{2}}|\nb|^{s+\frac{d-\alp}{2}} u }_{L^2}$ is equivalent with $\nrm{u}_{\dot{H}^{s+\frac{d-\alp}{2}}}$. To handle the general case of $\rho$ vanishing at infinity, we are forced to use the modified sum $\nrm{|\nb|^{s+\frac{d-\alp}{2}} u }_{L^2}^2 + \nrm{ \rho^{-\frac{1}{2}}|\nb|^s\rho}_{L^2}^2$. This requires the use of weighted Gagliardo--Nirenberg--Sobolev inequalities to control nonlinearity, which is given in Lemma \ref{lem:inequality}. 

\bigskip

\noindent \textbf{Local well-posedness with pressure} (Theorem \ref{thm:lwp-pressure}). When the system has pressure ($c_p>0$), at least on formal grounds the potential term is of lower order, which suggests local well-posedness regardless of the sign of $c_K$. We show rigorously that this is indeed the case. We emphasize that especially when $\alp$ approaches $d$, the potential term is of only ``slightly'' lower order than the other terms, and therefore cannot be handled by a simple energy estimate as in the Euler--Poisson case. The basic strategy of the local well-posedness proof is as follows. We first rewrite the system in terms of $u$ and $q = \rho^{\tilde{\gamma}}/\tilde{\gamma}$ with $\tilde{\gamma} = (\gmm-1)/2$ so that the system \eqref{eq:Euler-second} becomes  symmetric in the leading order. For the resulting equation, we observe a cancellation structure for the time derivative of the quantity \begin{equation*}
\begin{split}
\nrm{u}_{H^m}^2 + \nrm{q}_{H^m}^2 - Cc_K \nrm{ q^{\frac{1}{2\tilde{\gmm}}-1 } |\Lmb^{\frac{\alp-d}{2}}\rd^mq|}_{L^2}^2 .
\end{split}
\end{equation*} Inclusion of the last quantity is essential, to handle the loss of derivative coming from the potential term. Then, one can close an $H^m$ a priori estimate for the solution after bounding $\nrm{ q^{\frac{1}{2\tilde{\gmm}}-1 } |\Lmb^{\frac{\alp-d}{2}}\rd^mq|}_{L^2}$ in terms of $\nrm{q}_{H^m}$. Again, the sharp commutator estimate involving $\Lmb^{\eps}$ and a weighted Gagliardo--Nirenberg--Sobolev type inequality is crucial in this step. The unfortunate restriction $\gmm \le \frac{5}{3}$ stems from this procedure, since the commutator estimate requires certain smoothness of the coefficient, which is given by $q^{\frac{1}{2\tilde{\gmm}}-1}$.  One can expect that even in this case with pressure, local well-posedness proof simplify when either (i) $\rho$ attains a uniform positive lower bound or (ii) $c_K<0$. Indeed, when either of these assumptions is satisfied, the restriction $\gmm\le\frac{5}{3}$ can be dropped. This issue discussed in some more detail in Remark \ref{rem:pressure} below. 

\bigskip
 
\noindent \textbf{Finite-time singularity formation} (Theorems \ref{thm_att}, \ref{thm_rep}, and \ref{thm_iso}). The strategy for the formation of finite-time singularity relies on some energy estimates (cf. \cite{LPZ,S85,W14,Xi}). Note that the Euler--Riesz system \eqref{eq:Euler-second} conserves the mass, momentum, and total energy; see Lemma \ref{lem_en}. By employing these physical quantities, we investigate the time evolution of the moment of inertia or the internal energy for the Euler--Riesz system. When $\gamma > 1$ and $c_K > 0$ (with some assumption on the range of $\alp$), the total energy can be negative, which is preserved in time. We use this observation to show that under the hypothesis that a strong solution exists globally in time, the moment of inertia can be strictly negative in finite time, which is a contradiction. On the other hand, when $\gamma > 1$ and $c_K < 0$, the total energy is always positive, thus the previous argument cannot be applied to this case. In this case, roughly speaking, we estimate the lower and upper bound of the internal energy as follows:
\[
\frac{C_{l,1}}{(1+C_{l,2}t^2)^{\frac{d(\gamma-1)}{2}}} \leq \frac{1}{\gamma-1}\int_{\R^d} \rho^\gamma(t,x)\,dx \leq \frac{C_u}{(1+t)^{d(\gamma-1)}}
\]
for some positive constants $C_{l,1}$, $C_{l,2}$, and $C_u$ and all $t \geq 0$. This asserts that if $C_{l,1} > C_u C_{l,2}$, the life-span $T$ of the solution should be finite. In the isothermal pressure case, $\gamma=1$, the internal energy becomes $\int_{\R^d} \rho \ln \rho\,dx$, thus the total energy can be negative. However, we cannot use a similar argument as in the attractive case since $\rho \ln \rho$ is non-positive when $\rho \in [0,1]$. We show that the non-negative part of the internal energy can be controlled by the moment of inertia, and derive a second-order differential inequality for the moment of inertia. From the inequality, we obtain an explicit upper bound of the moment of inertia. This enables us to find some class of initial data leading to a formation of finite-time singularity. Our methods can be directly applied to the Coulomb case, $\alpha = d-2$, and we also would like to emphasize that the analysis of finite-time singularity formation for the multidimensional isothermal Euler--Poisson system with repulsive forces (without radially symmetry) has not been addressed so far, to the best of our knowledge.


\subsection{Organization of the paper}

The rest of this paper is organized as follows. In Section \ref{sec:lin}, we perform a formal linear analysis around a constant state to observe hyperbolic nature of the Euler--Poisson system, which suggests local well-posedness. Then in Sections \ref{sec:well1} and \ref{sec:well2}, we establish local well-posedness of classical solutions to the Euler--Riesz system without and with pressure, respectively. Finally, in Section \ref{sec:sing} we establish finite-time singularity formation for a class of large initial data.

\section{Linear analysis}\label{sec:lin}

As we mentioned earlier, a straightforward energy estimate does not close for the Euler--Riesz system \eqref{eq:Euler-second}.
To see whether this system has a chance to be well-posed, we take the linearization approach. Formally, $(\rho,u) = (1,0)$ defines a steady state to \eqref{eq:Euler-second}. Linearization around this steady state gives \begin{equation} \label{eq:Euler-second-linear}
\left\{
\begin{aligned}
&\rd_t\rho  = - \nabla \cdot u, \\
&\rd_t u = c_K\Lmb^{\alp-d} \nb \rho - c_p \nabla\rho . 
\end{aligned}
\right.
\end{equation} Assuming for simplicity that $\Omg = \bbT^d$ and taking the Fourier transform gives \begin{equation} \label{eq:Euler-second-linear-Fourier}
\left\{
\begin{aligned}
\frac{d}{dt}\widehat{\rho}(t,k) & = - ik \cdot \widehat{u}(t,k), \\
\frac{d}{dt}\widehat{u}(t,k) & = -ik(-c_KC_{\alp,d}|k|^{\alp-d} + c_p)\widehat{\rho}(t,k),
\end{aligned}
\right.
\end{equation} with some constant $C_{\alp,d}>0$. One sees from \eqref{eq:Euler-second-linear-Fourier} the following conservation law: \begin{equation*}
\begin{split}
\frac{1}{2} \frac{d}{dt} \left(   (-c_KC_{\alp,d} |k|^{\alp-d} + c_p) |\widehat{\rho}(t,k)|^2  + |\widehat{u}(t,k)|^2 \right) = 0. 
\end{split}
\end{equation*} This shows that 
\begin{itemize}
	\item If $c_p \ge 0$ and $c_K \le 0$, then the system \eqref{eq:Euler-second-linear} is well-posed in $H^s$-spaces for any $\alp$.
\end{itemize} We now consider the case $c_K > 0$. Then, \begin{equation*}
\begin{split}
\frac{d^2}{dt^2} \widehat{\rho}(t,k) = |k|^2\left( c_K C_{\alp,d} |k|^{\alp-d} -c_p \right)  \widehat{\rho}(t,k). 
\end{split}
\end{equation*} Since we assume that $\alp-d < 0$, for $|k|$ sufficiently large, we have that $c_K C_{\alp,d} |k|^{\alp-d} -c_p < 0$ as long as $c_p>0$. When $c_p = 0$, $|k|^{\alp-d+2} \rightarrow +\infty$ from $-2 < \alp -d$. Therefore, we conclude that \begin{itemize}
	\item If $c_K > 0$ and $c_p > 0$, the system \eqref{eq:Euler-second-linear} is still well-posed for $\alp-d < 0$.
	\item However, when $c_K > 0$ and $c_p = 0$, the system \eqref{eq:Euler-second-linear} is \textit{ill-posed} for $\alp-d>-2$. 
\end{itemize} 
In conclusion, assuming $-2<\alp-d<0$, linear analysis suggests local well-posedness when either $c_p>0$ \textit{regardless} of the sign of $c_K$ or $c_p = 0$ and $c_K \le 0$. Moreover, ill-posedness is expected when $c_p = 0$ and $c_K >0$. Similarly, even when $c_p > 0$, if $\nb_x\Lmb^{\alp-d}$ is replaced with $\Lmb^{\alp+1-d}$ (which is of the same order), ill-posedness should occur. In the following sections, we rigorously prove local well-posedness for the nonlinear system. We would like to emphasize that our proofs of well-posedness carry over (with minor modifications) to the case when the Riesz kernel \eqref{ker} is generalized to some kernel $K$ which satisfies the asymptotics $K(z) \simeq |z|^{-\alp}$ in the limit $|z|\rightarrow 0$, together with some differentiability properties.

%
%
%
%
%

\section{Well-posedness for the pressureless Euler--Riesz system}\label{sec:well1}

In this section, we assume that the domain is given by $\bbR^d$, where some delicate issues involving decay of $\rho$ at infinity appears. Modifying the proofs to the case of $\bbT^d$ is straightforward. Moreover, we always assume that $\rho$ is a strictly positive function in $\bbR^d$.


 We define the \textit{modified} $H^m$ ``norm" for $\rho$, as follows:   \begin{equation}\label{eq:modified-H-def}
\begin{split}
\nrm{\rho}_{\widetilde{H}^m}^2 := \sum_{0< |k|\le m} \nrm{\rho^{-\frac{1}{2}} \rd^k\rho }_{L^2}^2. 
\end{split}
\end{equation} Together with $\rho \in L^\infty$, the quantity $\widetilde{H}^m$ controls the usual Sobolev norm: \begin{equation*}
\begin{split}
\nrm{\nb\rho}_{H^{m-1}} \le \nrm{\rho}_{L^\infty}^{\frac{1}{2}}\nrm{\rho}_{\widetilde{H}^m}  .
\end{split}
\end{equation*} Since $\rho$ may decay at infinity, $\widetilde{H}^m$ encodes some extra decay of derivatives at infinity. Note that in \eqref{eq:modified-H-def}, the term with $k = 0$ is omitted. While we could have included this term in the definition of $\widetilde{H}^m$, doing so enforces that $\rho \in L^1(\bbR^d)$. On the other hand, our well-posedness result is able to treat the data for which $\rho$ converges to some positive constant at infinity. Once $d-2<\alp<d$ is fixed, we define the control quantity as \begin{equation*}
\begin{split}
\nrm{(u,\rho)}_{X^{m}} := \nrm{u}_{H^{m + \frac{d-\alp}{2}}} + \nrm{\rho }_{\widetilde{H}^m} + \nrm{ \rho  }_{ L^\infty } + \nrm{\nb(\ln\rho )}_{L^\infty}.
\end{split}
\end{equation*} We say that $(u,\rho)\in X^m$ if $\nrm{(u,\rho)}_{X^{m}}<+\infty$. 

\begin{theorem}[Local well-posedness of the pressureless and repulsive system]\label{thm:pressureless-lwp}
	{Assume that $d-2<\alp < d$.} Then the system \eqref{eq:Euler-second} in the pressureless and repulsive case is locally well-posed in the space $X^m$ for any $m> \frac{d}{2}+2$. More precisely, given any initial data satisfying \begin{equation*}
	\begin{split}
	\nrm{(u_0,\rho_0)}_{X^{m}}  = \nrm{u_0}_{ H^{m + \frac{d-\alp}{2}}} + \nrm{ \rho_0 }_{ L^\infty } + \nrm{\rho_0}_{\widetilde{H}^m} + \nrm{\nb(\ln\rho_0)}_{L^\infty} < +\infty,
	\end{split}
	\end{equation*} there exist $T >0$ and a unique solution $(u,\rho)$ to the system \eqref{eq:Euler-second} with $c_p = 0$ and $c_K < 0$ defined on $[0,T)$, satisfying the initial condition and $(u,\rho) \in L^\infty(0,T;X^{m})$; that is, \begin{equation*}
	\begin{split}
	 \nrm{u(t)}_{ H^{m + \frac{d-\alp}{2}}} + \nrm{ \rho(t) }_{ L^\infty } + \nrm{\rho(t)}_{\widetilde{H}^m} + \nrm{\nb(\ln\rho(t))}_{L^\infty} < +\infty,\quad t <T. 
	\end{split}
	\end{equation*}
\end{theorem}

We will need a few technical lemmas. 
\begin{lemma}[{{\cite{JK1}}}]\label{lem:inequality}
	Let $g >0$ on $\bbR^d$. Given an integer $m\ge 1$ and a $k$-tuple $(\ell_1, \cdots, \ell_k)$ of $d$-vectors satisfying $m = |\ell_1| + \cdots + |\ell_k|$, we have that \begin{equation}\label{eq:ineq-GNS-weight}
	\begin{split}
	\int_{\R^d} \frac{ 1 }{ g^{2k-1} } \prod_{1\le i \le k} |\rd^{\ell_i}g|^2 \,dx\lesssim_{m,d} \int_{\R^d} \frac{|\nb^mg|^2}{g}\,dx + \int_{\R^d} \frac{|\nb g|^{2m}}{g^{2m-1}}\,dx. 
	\end{split}
	\end{equation} In particular, we obtain that
	\begin{equation*}
	\begin{split}
	\nrm{ \nb^n ( \frac{\nb^{m-n}g}{\sqrt{g}} ) }_{L^2}^2 \lesssim_{m,d}  \int_{\R^d} \frac{|\nb^mg|^2}{g} \,dx+ \int_{\R^d} \frac{|\nb g|^{2m}}{g^{2m-1}}  \,dx
	\end{split}
	\end{equation*} for any integer $0\le n\le m$. 
\end{lemma}
 
\begin{proof}[Sketch of the proof]
	First, observe that it suffices to prove the inequality \begin{equation}\label{eq:ineq-GNS}
	\begin{split}
	\int_{\R^d} \frac{|\nb^\ell  g|^{\frac{2m}{\ell}}}{ g^{\frac{2m}{\ell}-1}}\,dx \lesssim_{m,d} \int_{\R^d} \frac{|\nb^m g|^2}{ g}\,dx + \int_{\R^d} \frac{|\nb g|^{2m}}{ g^{2m-1}}\,dx
	\end{split}
	\end{equation} for all integer $1 < \ell < m$. Once \eqref{eq:ineq-GNS} is proved, we then deduce \eqref{eq:ineq-GNS-weight} using H\"older's inequality: \begin{equation*}
	\begin{split}
	\int_{\R^d} \frac{ 1 }{  g^{2k-1} } \prod_{1\le i \le k} |\rd^{\ell_i} g|^2\,dx\le \prod_{1\le i\le k} \left(\int_{\R^d} \frac{|\nb^{|\ell_i|}  g|^{\frac{2m}{|\ell_i|}}}{ g^{\frac{2m}{|\ell_i|}-1}}\,dx\right)^{\frac{|\ell_i|}{m}} \lesssim_{m,d} \int_{\R^d} \frac{|\nb^m g|^2}{ g} \,dx+ \int_{\R^d} \frac{|\nb g|^{2m}}{ g^{2m-1}}\,dx.
	\end{split}
	\end{equation*} We show how to prove \eqref{eq:ineq-GNS} in the case $m =5$. Then, we need to prove that \begin{equation}\label{eq:goal}
	\begin{split}
	I_{\ell} \lesssim_{m,d} I_0 + I_5
	\end{split}
	\end{equation} for all $1<\ell<5$, where \begin{equation*}
	\begin{split}
	I_{\ell} := \int_{\R^d} \frac{|\nb^\ell  g|^{\frac{10}{\ell}}}{ g^{\frac{10}{\ell}-1}}\,dx.
	\end{split}
	\end{equation*} We fix some partial derivative $\rd^\ell$ of order $\ell$ and compute \begin{equation*}
	\begin{split}
	&\int_{\R^d} \frac{|\rd^\ell g|^{\frac{10}{\ell}-1}}{g^{^{\frac{10}{\ell}-1}}} \mathrm{sgn} (\rd^{\ell}g)  \rd^\ell g \,dx \cr
	&\quad = - ({\frac{10}{\ell}-1})\int_{\R^d} \frac{|\rd^\ell g|^{\frac{10}{\ell}-2}}{g^{\frac{10}{\ell}-1}} \, \rd^{\ell-1} g \, \rd^{\ell+1} g \,dx+ ({\frac{10}{\ell}-1} )\int_{\R^d} \frac{|\rd^\ell g|^{\frac{10}{\ell}-1}}{g^{\frac{10}{\ell}}} \, \rd g \, \rd^{\ell-1} g \, \mathrm{sgn}(\rd^\ell g)\,dx. 
	\end{split}
	\end{equation*} Note that $2.5\le \frac{10}{\ell} \le 5$. Using H\"older's inequality to the right hand side, we obtain that \begin{equation*}
	\begin{split}
	I_{\ell} \lesssim I_{\ell-1}^{ \frac{\ell-1}{10}} I_{\ell}^{1-\frac{\ell}{5}} I_{\ell+1}^{\frac{\ell+1}{10}} + I_1^{\frac{1}{10}} I_{\ell-1}^{\frac{\ell-1}{10}} I_{\ell}^{1-\frac{\ell}{10}}. 
	\end{split}
	\end{equation*}  With $\epsilon$-Young inequality, we obtain \begin{equation*}
	\begin{split}
	I_{\ell} \le \eps I_\ell + C_{\eps} ( I_1 + I_{\ell-1} + I_{\ell+1} ). 
	\end{split}
	\end{equation*} Combining this inequality in the cases $\ell = 2, 3, 4$, it is straightforward to obtain \eqref{eq:goal}. We omit the details. 
\end{proof}

\begin{lemma}\label{lem:comm-Lmb}
	Let $s\ge0$. For a vector field $v \in  (H^{\frac{d}{2}+1+s+\eps}(\bbR^d))^d$ and $f \in H^s(\bbR^d)$, we have the commutator estimate 
	\[
	\begin{split}
	\nrm{[\Lmb^s, v\cdot\nb]f}_{L^2} \lesssim_{s,d,\eps}\nrm{ v}_{H^{\frac{d}{2}+1+s+\eps}} \nrm{f}_{H^{s}}. 
	\end{split}
	\]
\end{lemma}
\begin{proof}
	
	We recall the following commutator estimate from \cite{CCCGW}:
	\begin{proposition}[{{\cite[Proposition 2.1]{CCCGW}}}]\label{prop:comm}
		For $s \in \bbR$ and any $\eps > 0$, there exists $C_{s,\eps}>0$ depending only on $s, \eps$ such that \begin{equation*}
		\begin{split}
		\nrm{[\Lmb^s \rd_{x_i}, g]f}_{L^2} \le C_{s,\eps}\left( \nrm{g}_{H^{\frac{d}{2}+1+\eps}} \nrm{\Lmb^s f}_{L^2} + \nrm{g}_{H^{\frac{d}{2}+1+s+\eps}} \nrm{f}_{L^2} \right)
		\end{split}
		\end{equation*} for any $i = 1, \cdots, d$.
	\end{proposition}
	Strictly speaking in \cite{CCCGW}, the proof is provided only in the case of $d = 2$, but the proof readily extends to any $d\ge 1$.  Now, observing\begin{equation*}
	\begin{split}
	[\Lmb^s \rd_{x_i}, g]f = [\Lmb^s , g\rd_{x_i}]f + \Lmb^s( f\rd_{x_i}g ),
	\end{split}
	\end{equation*} taking $g = v_i$ for each $i = 1, \cdots, d$, and summing up in $i$, we obtain the identity\begin{equation*}
	\begin{split}
	[\Lmb^s \nb, v]f = [\Lmb^s , v\cdot\nabla]f + \Lmb^s( f \nabla\cdot v ).
	\end{split}
	\end{equation*} Using the Fourier transform, it is not difficult to show that \begin{equation*}
	\begin{split}
	\nrm{\Lmb^s( f \nabla\cdot v )}_{L^2} \le C_{s,\eps}\left( \nrm{v}_{H^{\frac{d}{2}+1+\eps}} \nrm{\Lmb^s f}_{L^2} + \nrm{v}_{H^{\frac{d}{2}+1+s+\eps}} \nrm{f}_{L^2}  \right) 
	\end{split}
	\end{equation*} holds (cf. \cite[proof of Proposition 2.1]{CCCGW}).  Applying the above proposition, the estimate \begin{equation*}
	\begin{split}
	\nrm{[\Lmb^s,  v\cdot\nb ]f}_{L^2} \le C_{s,\eps}\left( \nrm{v}_{H^{\frac{d}{2}+1+\eps}} \nrm{\Lmb^s f}_{L^2} + \nrm{v}_{H^{\frac{d}{2}+1+s+\eps}} \nrm{f}_{L^2}  \right) \le C_{s,\eps} \nrm{v}_{H^{\frac{d}{2}+1+s+\eps}} \nrm{f}_{H^s}
	\end{split}
	\end{equation*} follows for $s\ge 0$. \end{proof}

\begin{proof}[Proof of Theorem \ref{thm:pressureless-lwp}]
 
In the proof, the constants may depend on $d,$ $m$, and $\alp$ but not on the solution. We divide the proof into a few steps.

\medskip

\noindent \textit{(i) a priori estimates}

\medskip

\noindent
Fix some derivative $\rd^m$ of order $m$ and we obtain \begin{equation*}
\begin{split}
\rd_t (\rd^mu) + u\cdot\nb \rd^m u + [\rd^m,u\cdot\nb] u = -\Lmb^{\alp-d} \nb\rd^m\rho.
\end{split}
\end{equation*} We now apply $\Lmb^{\frac{d-\alp}{2}}$ to both sides: \begin{equation*}
\begin{split}
\rd_t (\Lmb^{\frac{d-\alp}{2}}\rd^m u) + u\cdot\nb (\Lmb^{\frac{d-\alp}{2}}\rd^m u) + [\Lmb^{\frac{d-\alp}{2}}, u\cdot\nb]\rd^m u + \Lmb^{\frac{d-\alp}{2}}( [\rd^m,u\cdot\nb] u ) = -\Lmb^{\frac{\alp-d}{2}}\nb\rd^m\rho.
\end{split}
\end{equation*} Then, we have that\begin{equation*}
\begin{split}
\nrm{  [\Lmb^{\frac{d-\alp}{2}}, u\cdot\nb]\rd^m u }_{L^2} \le  C \nrm{ u }_{H^{m+\frac{d-\alp}{2}}}^2
\end{split}
\end{equation*} by Lemma \ref{lem:comm-Lmb}, since $m>\frac{d}{2}+2$. To handle the term $\Lmb^{\frac{d-\alp}{2}}( [\rd^m,u\cdot\nb] u )$, it suffices to estimate expressions of the form $\Lmb^{\frac{d-\alp}{2}}( \rd^k u\cdot\nb \rd^{m-k}u )$ where $1\le |k| \le m$ and $\rd^k$ and $\rd^{m-k}$ denote some derivatives of order $|k|$ and $m-|k|$, respectively. In the extreme cases $|k| = 1$ and $|k| = m$, we bound \begin{equation*}
\begin{split}
\nrm{ \Lmb^{\frac{d-\alp}{2}}( \rd u \cdot\nb \rd^{m-1}u )}_{L^2} +  
\nrm{ \Lmb^{\frac{d-\alp}{2}}( \rd^m u \cdot\nb u )}_{L^2} & \le C \nrm{ \rd u \cdot\nb \rd^{m-1}u }_{H^{\frac{d-\alp}{2}}} +  
\nrm{ \rd^m u \cdot\nb u }_{H^{\frac{d-\alp}{2}}} \\
& \le C\nrm{\widehat{\nb u} }_{L^1_\xi}\nrm{ u }_{H^{m+\frac{d-\alp}{2}}} + C\nrm{ |\xi|^{\frac{d-\alp}{2}} \widehat{\nb u}}_{L^1_\xi} \nrm{u}_{H^m} \\
& \le  C\nrm{ u }_{H^{m+\frac{d-\alp}{2}}}^2
\end{split}
\end{equation*} using the elementary product estimate in the Fourier space. For $1 < |k| < m$, one can simply estimate \begin{equation*}
\begin{split}
\nrm{ \Lmb^{\frac{d-\alp}{2}}( \rd^k u \cdot\nb \rd^{m-k}u )}_{L^2} \le C\nrm{  \rd^k u \cdot\nb \rd^{m-k}u }_{H^1} \le C\nrm{ u }_{H^{m+\frac{d-\alp}{2}}}^2
\end{split}
\end{equation*} using the Gagliardo--Nirenberg--Sobolev inequality. Here it was used that $m > \frac{d}{2} + 2$. Collecting the estimates, we obtain with $U_m := \Lmb^{\frac{d-\alp}{2}} \rd^m u$ and $R_m := \rd^m\rho$
that
\begin{equation}\label{eq:U_m}
\begin{split}
\left| \frac{1}{2} \frac{d}{dt} \nrm{U_m}_{L^2}^2 + \int_{\R^d} U_m \cdot  \Lmb^{\frac{\alp-d}{2}} \nb R_m \,dx \right| \le C \nrm{ u }_{H^{m+\frac{d-\alp}{2}}}^3. 
\end{split}
\end{equation} 
On the other hand, the equation for $R_m$ is given by \begin{equation*}
\begin{split}
\rd_tR_m + u\cdot\nb R_m + [\rd^m, u\cdot\nb]\rho = -\rho\nb\cdot \Lmb^{\frac{\alp-d}{2}}U_m - [\rd^m, \rho\nb\cdot] u . 
\end{split}
\end{equation*} Let us estimate \begin{equation*}
\begin{split}
 \frac{1}{2} \frac{d}{dt}\int_{\R^d} \frac{1}{\rho} R_m^2 \,dx& =  - \frac{1}{2}\int_{\R^d} \frac{\rd_t\rho}{\rho^2} R_m^2 \,dx + \int_{\R^d} \frac{1}{\rho} R_m\rd_t R_m \,dx=: I + J.
\end{split}
\end{equation*} First, $I$ can be easily estimated by \[
I \leq \frac12 \| \pa_t \ln \rho\|_{L^\infty}\int_{\R^d} \frac1\rho  R_m^2\,dx.
\]
Since 
\[
\pa_t \ln \rho = -u \cdot \nabla \ln \rho - \nabla \cdot u,
\]
we get
\[
I \leq \frac12 \lt(\|u\|_{L^\infty}\| \nabla \ln \rho\|_{L^\infty} + \|\nabla u\|_{L^\infty} \rt)\int_{\R^d} \frac1\rho  R_m^2\,dx \leq C\|u\|_{H^m}\lt(1 + \| \nabla \ln \rho\|_{L^\infty} \rt)\int_{\R^d} \frac1\rho  R_m^2\,dx.
\]
Next, we write $J = J_1 + J_2 + J_3 + J_4$ where \begin{equation*}
\begin{split}
J_1 = - \int_{\R^d} \frac{1}{\rho}R_m u\cdot\nb R_m \,dx, \quad
J_2 = -\int_{\R^d} \frac{1}{\rho} R_m ( [\rd^m,u\cdot\nb]\rho )\,dx,
\end{split}
\end{equation*}\begin{equation*}
\begin{split}
J_3 = - \int_{\R^d} R_m \nb\cdot \Lmb^{\frac{\alp-d}{2}}U_m\,dx ,
\quad \mbox{and} \quad 
J_4 = -\int_{\R^d} \frac1\rho R_m [\rd^m, \rho\nb\cdot] u \,dx.
\end{split} 
\end{equation*} We shall estimate $J_1, J_2,$ and $J_4$. To begin with, \begin{equation*}
\begin{split}
J_1  &= \frac12\int_{\R^d} \nabla \cdot \lt(\frac{u}{\rho} \rt) R_m^2\,dx = \frac12\int_{\R^d} (\nabla \cdot u - u \nabla \ln \rho) \frac1\rho  R_m^2\,dx\cr
&\leq C\|u\|_{H^m}\lt(1 + \| \nabla \ln \rho\|_{L^\infty} \rt)\int_{\R^d} \frac1\rho  R_m^2\,dx. 
\end{split}
\end{equation*} To handle $J_2$,  we need to estimate \begin{equation*}
\begin{split}
\int_{\R^d} \frac{1}{\rho} (\rd^m\rho) \rd^\ell u \cdot \nb \rd^{m-\ell}\rho\,dx , \quad 0 < \ell \le m.
\end{split}
\end{equation*} We consider two cases: if $\ell < m- \frac{d}{2}$, then we simply estimate \begin{equation*}
\begin{split}
\left| \int_{\R^d} \frac{1}{\rho} (\rd^m\rho) \rd^\ell u \cdot \nb \rd^{m-\ell}\rho  \,dx\right| \le C\nrm{\rd^{\ell}u}_{L^\infty} \nrm{  \frac{\rd^m\rho}{\sqrt{\rho}} }_{L^2} \nrm{  \frac{\rd^{m-\ell+1}\rho}{\sqrt{\rho}} }_{L^2} \le C \nrm{u}_{H^m} \nrm{  \frac{\rd^m\rho}{\sqrt{\rho}} }_{L^2} \nrm{  \frac{\rd^{m-\ell+1}\rho}{\sqrt{\rho}} }_{L^2} 
\end{split}
\end{equation*} with Sobolev embedding. Next, when $\ell \ge m- \frac{d}{2}$, we define $s$ to be the smallest integer satisfying $s \ge \frac{d}{2} + \ell -m$. Setting $1 < p, q <\infty$ to satisfy \begin{equation*}
\begin{split}
m-\ell = \frac{d}{2} (1 - \frac{1}{p}) = \frac{d}{2q}, 
\end{split}
\end{equation*} we have first by H\"older inequality that \begin{equation*}
\begin{split}
\left| \int_{\R^d} \frac{1}{\rho} (\rd^m\rho) \rd^\ell u \cdot \nb \rd^{m-\ell}\rho \,dx \right| \le C\nrm{  \frac{\rd^m\rho}{\sqrt{\rho}} }_{L^2}\nrm{\rd^\ell u}_{L^{2p}} \nrm{ \frac{1}{\sqrt{\rho}} \rd^{m+1-\ell}\rho }_{L^{2q}},
\end{split}
\end{equation*} and then using Gagliardo-Nirenberg-Sobolev inequality and Lemma \ref{lem:inequality}, we have \begin{equation*}
\begin{split}
\nrm{\rd^\ell u}_{L^{2p}} \nrm{ \frac{1}{\sqrt{\rho}} \rd^{m+1-\ell}\rho }_{L^{2q}} &\le C \nrm{u}_{H^m}\left( \sum_{i=0}^{s} \nrm{ \rd^i  ( \frac{1}{\sqrt{\rho}} \rd^{m+1-\ell}\rho  )  }_{L^2}  \right)\\
&\le C \nrm{u}_{H^m}\left( \sum_{j=1}^{s + m + 1 - \ell} \left( \int_{\R^d} \frac{|\nb^j \rho|^2}{\rho} \,dx+ \int_{\R^d} \frac{|\nb \rho|^{2j}}{\rho^{2j-1}} \,dx\right)^{\frac{1}{2}}  \right).
\end{split}
\end{equation*} Note that $s + m + 1 - \ell < \frac{d}{2} + 2 $. Moreover, we can bound \begin{equation*}
\begin{split}
\int_{\R^d} \frac{|\nb \rho|^{2j}}{\rho^{2j-1}} \,dx\le C \nrm{ \nb (\ln\rho)}_{L^\infty}^{2(j-1)} \int_{\R^d} \frac{|\nb\rho|^2}{\rho}\,dx.
\end{split} 
\end{equation*} Therefore, we obtain that \begin{equation*}
\begin{split}
|J_2| &\le C \nrm{u}_{H^m} \left( \sum_{0<k\le m} (   \nrm{\frac{1}{\sqrt{\rho} }R_k }_{L^2}^2    + \nrm{\nb(\ln\rho)}_{L^\infty}^{2(k-1)}  \nrm{\frac{1}{\sqrt{\rho} }R_1 }_{L^2}^2 )  \right)\\
&\le C (1 + \nrm{\nb(\ln\rho)}_{L^\infty} )^{2(m-1)} \nrm{u}_{H^m}  \sum_{0<k\le m}  \nrm{\frac{1}{\sqrt{\rho} }R_k }_{L^2}^2 . 
\end{split}
\end{equation*} Now, we observe that to estimate $J_4$, it suffices to treat terms of the form \begin{equation*}
\begin{split}
\int_{\R^d} \frac{1}{\rho} R_m \rd^{\ell}\rho \rd^{m-\ell} \nb\cdot u\,dx
\end{split}
\end{equation*} for $\ell>0$, but it can be estimated exactly the same way with $J_2$. We conclude that \begin{equation}\label{eq:R_m}
\begin{split}
\left| \frac{1}{2} \frac{d}{dt}\int_{\R^d} \frac{1}{\rho} R_m^2 \,dx+ \int_{\R^d} R_m \nb\cdot \Lmb^{\frac{\alp-d}{2}}U_m \,dx \right| \le C \|u\|_{H^m}\lt(1 + \| \nabla \ln \rho\|_{L^\infty} \rt)^{2(m-1)}  \sum_{0<k\le m}  \nrm{\frac{1}{\sqrt{\rho} }R_k }_{L^2}^2  .
\end{split}
\end{equation} Using \eqref{eq:R_m} together with \eqref{eq:U_m}, we obtain with \begin{equation*}
\begin{split}
- \int_{\R^d} U_m \cdot \Lmb^{\frac{\alp-d}{2}} \nb R_m\,dx - \int_{\R^d} R_m \nb\cdot \Lmb^{\frac{\alp-d}{2}}U_m \,dx= 0 
\end{split}
\end{equation*} that \begin{equation*}
\begin{split}
\left| \frac{1}{2} \frac{d}{dt}\int_{\R^d}  U_m^2 \,dx+ \frac{1}{\rho} R_m^2 \,dx\right| \le C\lt(1 + \| \nabla \ln \rho\|_{L^\infty} \rt)^{2(m-1)}\nrm{ u }_{H^{m+\frac{d-\alp}{2}}} \left(\nrm{ u }_{H^{m+\frac{d-\alp}{2}}}^2 +  \sum_{0<k\le m}  \nrm{\frac{1}{\sqrt{\rho} }R_k }_{L^2}^2  \right).
\end{split}
\end{equation*} Repeating the argument for all possible partial derivatives of order $\le m$, we obtain \begin{equation}\label{eq:pressureless-apriori1}
\begin{split}
\left| \frac{d}{dt} \left( \nrm{u}^2_{H^{m+\frac{d-\alp}{2}}} + \nrm{\rho}_{\widetilde{H}^m}^2 \right) \right| \le C \lt(1 + \| \nabla \ln \rho\|_{L^\infty} \rt)^{2(m-1)}\nrm{ u }_{H^{m+\frac{d-\alp}{2}}} \left( \nrm{u}^2_{H^{m+\frac{d-\alp}{2}}} + \nrm{\rho}^2_{\widetilde{H}^m} \right) . 
\end{split}
\end{equation} Now, from the equation for $\rho$ \begin{equation*}
\begin{split}
\rd_t\rho + u\cdot\nb\rho = -\rho \nabla\cdot u ,
\end{split}
\end{equation*} we see that \begin{equation*}
\begin{split}
\rd_t \nb (\ln\rho) + u\cdot\nb \nb (\ln\rho) = - (\nb u)^T \nb (\ln\rho) - \nb (\nb\cdot u). 
\end{split}
\end{equation*} Evaluating the previous equation along the flow generated by $u$, we see that \begin{equation}\label{eq:pressureless-apriori2}
\begin{split}
\frac{d}{dt} \nrm{\nb (\ln\rho)}_{L^\infty} \le C\left( \nrm{\nb u}_{L^\infty}\nrm{\nb (\ln\rho)}_{L^\infty} + \nrm{\nb^2 u}_{L^\infty} \right) \le C \nrm{u}_{H^{m+\frac{d-\alp}{2}}} (1 + \nrm{\nb (\ln\rho)}_{L^\infty}) . 
\end{split}
\end{equation}   On the other hand, from the equation for $\rho$, we have \begin{equation}\label{eq:pressureless-apriori4}
\begin{split}
\frac{d}{dt} \nrm{\rho}_{L^\infty} \le C \nrm{\nb u}_{L^\infty} \nrm{\rho}_{L^\infty} \le C \nrm{ u }_{H^{m+\frac{d-\alp}{2}}}  \nrm{\rho}_{L^\infty}. 
\end{split}
\end{equation} We now recall the definition of the $X^m$-norm: \begin{equation*}
\begin{split}
\nrm{(u(t),\rho(t))}_{X^m} = \nrm{u(t)}_{ H^{m + \frac{d-\alp}{2}}} + \nrm{ \rho(t) }_{ L^\infty } + \nrm{\rho(t)}_{\widetilde{H}^m} + \nrm{\nb(\ln\rho(t))}_{L^\infty}.
\end{split}
\end{equation*} Combining \eqref{eq:pressureless-apriori1}, \eqref{eq:pressureless-apriori2}, and \eqref{eq:pressureless-apriori4}, we obtain that \begin{equation*}
\begin{split}
\left|\frac{d}{dt} \nrm{(u(t),\rho(t))}_{X^m}^2  \right|\le  C(1 + \nrm{(u(t),\rho(t))}_{X^m})^{2m+1}.
\end{split}
\end{equation*} This gives an a priori estimate for the solution of \eqref{eq:Euler-second} in the pressureless and repulsive case. 

\medskip

\noindent \textit{(ii) uniqueness}

\medskip

\noindent Assume that for some interval of time $[0,T]$, there exist two solutions $(\rho_1,u_1)$ and $(\rho_2,u_2)$ of \eqref{eq:Euler-second} with the same initial data $(\rho_0,u_0)\in X^m$, satisfying \begin{equation*}
\begin{split}
\max_{i=1,2}\sup_{t\in[0,T]}(\nrm{u_i(t)}_{ H^{m + \frac{d-\alp}{2}}} + \nrm{ \rho_i(t) }_{ L^\infty } + \nrm{\rho_i(t)}_{\widetilde{H}^m} + \nrm{\nb(\ln\rho_i(t))}_{L^\infty} ) = M < +\infty. 
\end{split}
\end{equation*} In the estimates below, the constant $C$ may depend on $M$ as well.  We define $\trho = \rho_1 - \rho_2$ and $\tu = u_1 - u_2$. Then, the equations for $\trho$ and $\tu$ read
\begin{equation}\label{eq:trho}
\begin{split}
\rd_t\trho + u_1\cdot\nb \trho + \tu \cdot\nb \rho_2 = -\trho \nb\cdot u_1 - {\rho_2\nb\cdot\tu }
\end{split}
\end{equation} and \begin{equation}\label{eq:tu}
\begin{split}
\rd_t \tu + u_1\cdot\nb \tu + \tu \cdot\nb u_2 = -\nb \Lmb^{\alp-d} \trho. 
\end{split}
\end{equation} Before we proceed, let us observe that the ratio $\rho_1/\rho_2$ remains bounded from above and below. To prove this, we simply compute using the equations for $\rho_1$ and $\rho_2$ that \begin{equation*}
\begin{split}
\frac{d}{dt} \left( \frac{\rho_1}{\rho_2} \right) = - {\left( \frac{u_1\cdot\nb\rho_1}{\rho_1} + \nabla\cdot u_1 + \frac{\rd_t\rho_2}{\rho_2}  \right)\left( \frac{\rho_1}{\rho_2} \right)}.
\end{split}
\end{equation*} Since we have a uniform pointwise estimate \begin{equation*}
\begin{split}
\left|  \frac{u_1\cdot\nb\rho_1}{\rho_1} + \nabla\cdot u_1 + \frac{\rd_t\rho_2}{\rho_2}  \right|(t,x) \le C, 
\end{split}
\end{equation*} we obtain that \begin{equation*}
\begin{split}
\exp(-Ct) \le\left( \frac{\rho_1}{\rho_2} \right)(t,x) \le \exp(Ct)
\end{split}
\end{equation*} uniformly for all $x \in \bbR^d$, recalling that $\rho_1 = \rho_2$ at $t = 0$. 
Assuming for simplicity that $\rho_0\in L^1(\bbR^d)$\footnote{To drop this assumption, one could either perform a suitably weighted $L^2$ estimate for $\trho$ or carry out higher norm estimate for $(\trho,\tu)$.} (this property propagates in time), we obtain in particular that, the quantity $\rho_2^{-1}|\trho|^2$ is integrable. For simplicity, we set $\eps = (d-\alp)/2 >0$. From \eqref{eq:trho}, we compute that 
\bq\label{eq:trho2}
\begin{split}
\frac{1}{2} \frac{d}{dt} \int_{\R^d} \frac{1}{\rho_2}|\trho|^2 \,dx &= \frac{1}{2}\int_{\R^d} \frac{\rd_t\rho_2}{\rho_2} \frac{1}{\rho_2}|\trho|^2 \,dx-\frac{1}{2} \int_{\R^d} \nb\cdot u_1 \frac{1}{\rho_2}|\trho|^2\,dx -\frac12 \int_{\R^d} \frac{u_1 \cdot \nb \rho_2}{\rho_2} \frac1{\rho_2} |\trho|^2\,dx\cr
&\quad + \int_{\R^d} \frac{\nb\rho_2}{\rho_2} \cdot \tu \trho\,dx- {\int_{\R^d} \trho\nb \cdot \tu}\,dx .
\end{split}
\eq
Similarly, using \eqref{eq:tu}, we compute 
\begin{equation}\label{eq:tu2}
\begin{split}
\frac{1}{2} \frac{d}{dt} \int_{\R^d}  |\Lmb^\eps \tu|^2 \,dx&=  - \int_{\R^d} \frac{1}{2} \nb\cdot u_1  |\Lmb^\eps \tu|^2 \,dx- \int_{\R^d}  \Lmb^\eps \tu \cdot [\Lmb^\eps,u_1\cdot\nb] \tu \,dx \cr
&\quad -  {\int_{\R^d}  \Lmb^\eps \tu \cdot \Lmb^\eps(\tu\cdot\nb u_2)} \,dx- \int_{\R^d} \Lmb^\eps \tu \cdot \nb \Lmb^{-\eps} \trho\,dx. 
\end{split}
\end{equation} We rewrite the last terms on the right hand sides of \eqref{eq:trho2} and \eqref{eq:tu2} as follows:\begin{equation*}
\begin{split}
- \int_{\R^d} \nb \cdot \tu  \trho \,dx =  \int_{\R^d}  \tu\cdot\nb\trho \,dx \quad \mbox{and} \quad - \int_{\R^d}  \Lmb^\eps \tu \cdot \nb \Lmb^{-\eps} \trho\,dx = - \int_{\R^d} \tu \cdot  \nb\trho \,dx . 
\end{split}
\end{equation*}  Then, we obtain \begin{equation}\label{eq:ineq1}
\begin{split}
\left| \frac{1}{2} \frac{d}{dt} \int_{\R^d} \frac{1}{\rho_2}|\trho|^2\,dx - \int_{\R^d} \tu\cdot\nb\trho \,dx  \right| \le C \nrm{\rho_2^{-\frac{1}{2}}\trho}_{L^2} ( \nrm{\rho_2^{-\frac{1}{2}}\trho}_{L^2} + \nrm{\tu}_{L^2})
\end{split}
\end{equation} and estimating the second to last term in \eqref{eq:tu2} by \begin{equation*}
\begin{split}
\lt|\int_{\R^d}  \Lmb^\eps \tu \cdot \Lmb^\eps(\tu\cdot\nb u_2)\,dx\rt| \leq \|\Lmb^\eps \tu\|_{L^2} \|\Lmb^\eps(\tu\cdot\nb u_2)\|_{L^2} \leq C\|\Lmb^\eps \tu\|_{L^2}\lt(\|\Lmb^\eps \tu\|_{L^2} + \|\tu\|_{L^2}\rt)
\end{split}
\end{equation*} (using Lemma \ref{lem:comm-Lmb}), we obtain \begin{equation}\label{eq:ineq2}
\begin{split}
\left| \frac{1}{2} \frac{d}{dt} \int_{\R^d}  |\Lmb^\eps \tu|^2\,dx  + \int_{\R^d} \tu\cdot\nb\trho \,dx \right| \le C\nrm{\Lmb^\eps \tu}_{L^2}^2.
\end{split}
\end{equation}   
Adding the inequalities \eqref{eq:ineq1} and \eqref{eq:ineq2}, we obtain \begin{equation*}
\begin{split}
\left| \frac{d}{dt} ( \nrm{\rho_2^{-\frac{1}{2}}\trho}_{L^2} + \nrm{\Lmb^\eps\tu}_{L^2} ) \right| \le C (\nrm{\rho_2^{-\frac{1}{2}}\trho}_{L^2}  + \nrm{\tu}_{L^2} + \nrm{\Lmb^\eps\tu}_{L^2}  ).
\end{split}
\end{equation*} Proceeding similarly for the equations for $\nb\trho$ and $\Lmb^\eps \nb\tu$, we can obtain \begin{equation*}
\begin{split}
\left| \frac{d}{dt} ( \nrm{\rho_2^{-\frac{1}{2}}\nb\trho}_{L^2} + \nrm{\Lmb^\eps\nb\tu}_{L^2} ) \right| \le C (\nrm{\rho_2^{-\frac{1}{2}}\trho}_{L^2}  + \nrm{\rho_2^{-\frac{1}{2}}\nb\trho}_{L^2}  + \nrm{\tu}_{L^2} + \nrm{\Lmb^\eps\tu}_{L^2}  + \nrm{\Lmb^\eps\nb\tu}_{L^2}  ).
\end{split}
\end{equation*} Finally, from the equation for $\tu$, it follows that \begin{equation*}
\begin{split}
\left| \frac{d}{dt} \nrm{\tu}_{L^2} \right| \le C(\nrm{\tu}_{L^2} + \nrm{\Lmb^{\alp-d}\nb\trho}_{L^2}  ) \le C(\nrm{\tu}_{L^2} + \nrm{\trho}_{H^1}  ) \le C(\nrm{\tu}_{L^2} + \nrm{\rho_2^{-\frac{1}{2}}\trho}_{L^2}  + \nrm{\rho_2^{-\frac{1}{2}}\nb\trho}_{L^2}  ).
\end{split}
\end{equation*} 
Here,  $\nrm{\Lmb^{\alp-d}\nb\trho}_{L^2}  \leq C\nrm{\trho}_{H^1}$ follows from $\alp-d<0$. Hence, for the quantity \begin{equation*}
\begin{split}
X = \nrm{\rho_2^{-\frac{1}{2}}\trho}_{L^2}  + \nrm{\rho_2^{-\frac{1}{2}}\nb\trho}_{L^2}  + \nrm{\tu}_{L^2} + \nrm{\Lmb^\eps\tu}_{L^2}  + \nrm{\Lmb^\eps\nb\tu}_{L^2}  ,
\end{split}
\end{equation*} we obtain the differential inequality \begin{equation*}
\begin{split}
\left| \frac{d}{dt} X \right| \le CX, 
\end{split}
\end{equation*} which guarantees that if $X(t=0) = 0$, then $X = 0$ for $t\in[0,T]$. This shows that $u_1 = u_2$ and $\rho_1 = \rho_2$. 

\medskip

\noindent \textit{(iii) existence}

\medskip

\noindent It only remains to prove existence of a solution to \eqref{eq:Euler-second} satisfying the above a priori estimates. For this purpose we fix some small time interval $[0,T_1]$ on which we have $X_m(t) \le 2X_m(0)$. Once  existence is shown in $[0,T_1]$, one can extend the time interval for existence as long as the a priori estimate does not blow up.

To show existence, we fix some initial data $(u_0,\rho_0)$ satisfying the assumptions of Theorem \ref{thm:pressureless-lwp} and consider the following \textit{viscous} system: given a parameter $\epsilon>0$, we consider \begin{equation}\label{eq:existence-viscous}
\begin{split}
&\rd_t \rho^{(\eps)} + u^{(\eps)}\cdot \nb \rho^{(\eps)} = - \rho^{(\eps)} \nabla\cdot u^{(\eps)} + \eps \lap \rho^{(\eps)},\\
&\rd_t u^{(\eps)} + u^{(\eps)}\cdot\nb u^{(\eps)} = - \nb \Lmb^{\alp-d} \rho^{(\eps)} + \eps \lap u^{(\eps)} 
\end{split}
\end{equation} with $C^\infty$ initial data \begin{equation}\label{eq:existence-initial}
\begin{split}
u_0^{(\eps)} = \varphi_\eps * u_0, \quad \rho_0^{(\eps)} = \varphi_\eps * \rho_0.
\end{split}
\end{equation} Here, $\varphi_\eps(x) = \eps^{-d}\varphi(\eps^{-1}x)$ is a smooth approximation of the identity. Existence of a local in time smooth solution to the system \eqref{eq:existence-viscous}--\eqref{eq:existence-initial} follows from a standard contraction mapping argument (\cite{Ka}), using the mild formulation and estimates for the heat semigroup. To be more precise, one can rewrite \eqref{eq:existence-viscous} as 
\begin{equation*}
\begin{split}
\rho^{(\eps)}(t) &= e^{\eps t\lap} \rho_0^{(\eps)} - \int_0^t e^{\eps(t-s)\lap} \nb \cdot ( \rho^{(\eps)} u^{(\eps)} ) \,ds ,\\
u^{(\eps)}(t) &= e^{\eps t\lap} u_0^{(\eps)} - \int_0^t e^{\eps(t-s)\lap} (u^{(\eps)}\cdot\nb u^{(\eps)} + \nb \Lmb^{\alp-d} \rho^{(\eps)} ) \, ds 
\end{split}
\end{equation*} and prove that there exist $T = T(\eps, d, \rho_0, u_0)>0$ such that the operator \begin{equation*}
\begin{split}
(\rho^{(\eps)}, u^{(\eps)}) \mapsto \left( e^{\eps t\lap} \rho_0^{(\eps)} - \int_0^t e^{\eps(t-s)\lap} \nb \cdot ( \rho^{(\eps)} u^{(\eps)} ) \,ds, e^{\eps t\lap} u_0^{(\eps)} - \int_0^t e^{\eps(t-s)\lap} (u^{(\eps)}\cdot\nb u^{(\eps)} + \nb \Lmb^{\alp-d} \rho^{(\eps)} ) \, ds  \right)
\end{split}
\end{equation*} is a contraction mapping on the space $L^\infty([0,T]; B_{H^1}(2(\nrm{u_0}_{H^1} + \nrm{\rho_0}_{H^1})) )$. Here, $B_{H^1}(R)$ denotes the open ball of radius $R$ in the space $H^1(\bbR^d)$.

While a priori the time of existence becomes small as $\eps\to 0$, one can show that the a priori estimates for the inviscid equation proved above are satisfied for the viscous solutions as well. This guarantees uniform time of existence for the viscous solutions $(\rho^{(\eps)},u^{(\eps)})$. To see this for the case of $\nrm{ \nb (\ln \rho^{(\eps)}) }_{L^\infty}$, we consider the equation for $\ln\rho^{(\eps)}$: \begin{equation*}
\begin{split}
\rd_t \ln\rho^{(\eps)} + u^{(\eps)} \cdot \nb \ln\rho^{(\eps)} = -\nb\cdot u^{(\eps)} + \eps\left( \lap \ln\rho^{(\eps)} + |\nb \ln\rho^{(\eps)}|^2 \right).
\end{split}
\end{equation*} Taking a partial derivative $\rd_{x_i}$ for $i \in \{1,\cdots, d\}$, we find
\begin{equation*}
\begin{split}
\rd_t \rd_{x_i}\ln\rho^{(\eps)} + u^{(\eps)} \cdot \nb \rd_{x_i}\ln\rho^{(\eps)} &= -\nb\cdot \rd_{x_i}u^{(\eps)} - \rd_{x_i}u^{(\eps)} \cdot \nb\ln\rho^{(\eps)} \cr
&\quad + \eps\left( \lap \rd_{x_i}\ln\rho^{(\eps)} + 2 \nb \rd_{x_i}\ln\rho^{(\eps)}  \cdot \nb \ln\rho^{(\eps)}   \right).
\end{split}
\end{equation*} 
Given some $t$, at any local maximum point of $\rd_{x_i}\ln\rho^{(\eps)}$, we have \begin{equation*}
\begin{split}
\lap \rd_{x_i}\ln\rho^{(\eps)} \le 0, \quad  \nb \rd_{x_i}\ln\rho^{(\eps)}  = 0. 
\end{split}
\end{equation*} Arguing similarly for local minima of $\rd_{x_i}\ln\rho^{(\eps)}$, we can deduce the estimate \begin{equation*}
\begin{split}
\frac{d}{dt} \nrm{ \rd_{x_i}\ln\rho^{(\eps)} }_{L^\infty} \le C \nrm{\nb^2u^{(\eps)}}_{L^{\infty}} + C\nrm{\nb u^{(\eps)}}_{L^\infty} \nrm{ \nb\ln\rho^{(\eps)} }_{L^\infty} 
\end{split}
\end{equation*} for each $i = 1, \cdots, d$. This gives the a priori estimate for $\nrm{ \nb (\ln \rho^{(\eps)}) }_{L^\infty}$ which is uniform in the limit $\eps\to 0$. Repeating a similar argument for the other quantities, one obtains a sequence of solutions $( u^{(\eps)}, \rho^{(\eps)} )$ defined on $[0,T_1]$, with \begin{equation*}
\begin{split}
\nrm{(u^{(\eps)}, \rho^{(\eps)} )}_{L^\infty(0,T_1;X^m)} = \sup_{t \in [0,T_1] }\left( \nrm{u^{(\eps)}(t)}_{ H^{m + \frac{d-\alp}{2}}} + \nrm{ \rho^{(\eps)}(t) }_{ L^\infty } + \nrm{\rho^{(\eps)}(t)}_{\widetilde{H}^m} + \nrm{\nb(\ln\rho^{(\eps)}(t))}_{L^\infty} \right)
\end{split}
\end{equation*} uniformly bounded in $\eps$. Therefore, one can pass to a weakly convergent subsequence as $\eps\to 0$, with some limit $(u, \rho)$ satisfying \begin{equation*}
\begin{split}
\sup_{t \in [0,T_1] }\left( \nrm{u(t)}_{ H^{m + \frac{d-\alp}{2}}} + \nrm{ \rho(t) }_{ L^\infty } + \nrm{\rho(t)}_{\widetilde{H}^m} + \nrm{\nb(\ln\rho(t))}_{L^\infty} \right) < \infty.
\end{split}
\end{equation*} It is not difficult to show that $(u,\rho)(t=0) = (u_0,\rho_0)$ and that $(u,\rho)$ is a solution to \eqref{eq:Euler-second}. We omit the details. 
\end{proof}

%
%
%
%
%
\section{Well-posedness for the Euler--Riesz system}\label{sec:well2}

In this section, we consider the system \eqref{eq:Euler-second} with pressure; given some $\gmm\ge 1$, we fix without loss of generality that $c_p = 1/\gmm$. In the case $\gmm>1$, by introducing the variable 
\begin{equation*}
\begin{split}
q = \frac{1}{\tilde{\gmm}} \rho^{\tilde{\gmm}}, \quad \tilde{\gmm} = \frac{\gmm-1}{2}, 
\end{split}
\end{equation*} we have that \eqref{eq:Euler-second} turns into \begin{equation}\label{eq:Euler-Riesz-pressure2}
\left\{
\begin{aligned} 
&\rd_t q + u \cdot \nb q = - \tilde{\gmm} q \nb\cdot u, \\
&\rd_t u + u\cdot\nb u = - \tilde{\gmm} q \nb q + \tilde{c}_K \nb \Lmb^{\alp-d} ( q^{\frac{1}{\tilde{\gmm}}} )
\end{aligned}
\right.
\end{equation} with $\tilde{c}_K = c_K (\tilde{\gmm})^{\frac{1}{\tilde{\gmm}}}.$ Let us now consider the case $\gamma = 1$, which corresponds to the isothermal pressure law. In this case, we reformulate the system \eqref{eq:Euler-second} by introducing $q = \ln \rho$ (and $\tilde{c}_K = c_K$) as 
\begin{equation}\label{eq:Euler-Riesz-pressure3}
\left\{
\begin{aligned} 
&\rd_t q + u \cdot \nb q = -  \nb\cdot u, \\
&\rd_t u + u\cdot\nb u = - \nb q + \tilde{c}_K \nb \Lmb^{\alp-d} ( e^q ).
\end{aligned}
\right.
\end{equation} 

\begin{theorem}\label{thm:lwp-pressure}
	Assume that $1 < \gmm \le \frac{5}{3}$. For any $\tilde{c}_K\in\bbR$, the system \eqref{eq:Euler-Riesz-pressure2} is locally well-posed; for any initial data $(u_0,q_0)$ belonging to $H^m(\bbR^d)$ with $m>\frac{d}{2}+1$ and $\nb (\ln q_0) \in L^\infty(\bbR^d)$, there exist $T>0$ and a unique solution $(u,q)$ to \eqref{eq:Euler-Riesz-pressure2} satisfying the initial condition and the bounds \begin{equation*}
	\begin{split}
	\nrm{u(t)}_{H^m} + \nrm{q(t)}_{H^m} +\nrm{\nb (\ln q(t))}_{L^\infty} < +\infty
	\end{split}
	\end{equation*}  for $t<T$. In the case $\gmm = 1$, the system \eqref{eq:Euler-Riesz-pressure3} is locally well-posed with $(u,q)$ belonging to $H^m(\bbR^d)$ ($m>\frac{d}{2}+1$). That is, the assumption $\nb (\ln q_0) \in L^\infty(\bbR^d)$ is not necessary in this case. 
\end{theorem}

We need a lemma whose proof is completely parallel to that of Lemma \ref{lem:inequality}. 
\begin{lemma}[{{\cite{JK1}}}]\label{lem:power-Sobolev}
	Let $g>0$ on $\bbR^d$ and $\beta>0$. Then for any integer $k\ge 1$, \begin{equation*}
	\begin{split}
	\nrm{g^{\beta}}_{\dot{H}^{k}}^2 \lesssim_{\beta,d,k} \nrm{g^{\beta-1}|\nb^kg|}_{L^2}^2 + \nrm{ g^{-(1-\frac{\beta}{k})} |\nb g| }_{L^{2k}}. 
	\end{split}
	\end{equation*} In particular, for $g$ satisfying $|g^{-1}\nb g| \lesssim 1$ and $g^\beta \in L^2$, we have \begin{equation*}
	\begin{split}
	\nrm{g^{\beta}}_{\dot{H}^{k}}^2 \lesssim_{\beta,d,k} \nrm{g^{\beta-1}|\nb^kg|}_{L^2}^2 + \nrm{\nb (\ln g)}_{L^\infty}^{2k} \nrm{g^{\beta}}_{L^2}^{2}. 
	\end{split}
	\end{equation*}
\end{lemma}

\begin{proof}[Proof of Theorem \ref{thm:lwp-pressure}]
We first consider the case $\gmm>1$. Let us establish a priori estimates for a solution $(u,q)$ to \eqref{eq:Euler-Riesz-pressure2}. Straightforward computation yields
\begin{equation*}
\begin{split}
&\frac{d}{dt}\left( \frac{1}{2}\int_{\R^d} |\rd^mu|^2 + |\rd^m q|^2\,dx \right)- \tilde{c}_K \int_{\R^d} \rd^m u\cdot \rd^m \nb \Lmb^{\alp-d}(q^{\frac{1}{\tilde\gmm}})  \,dx \cr
&\quad = \frac12\int_{\R^d} (\nabla \cdot u) (|\pa^m q|^2 + |\pa^m u|^2)\,dx + \tilde \gamma \int_{\R^d} (\pa^m u) \nabla q \cdot \pa^m q\,dx\cr
&\qquad - \int_{\R^d} \pa^m q \cdot [\pa^m, u \cdot \nabla] q\,dx - \int_{\R^d} \pa^m u \cdot[\pa^m, u \cdot \nabla] u\,dx \cr
&\qquad -  \tilde \gamma\int_{\R^d} \pa^m q \cdot [\pa^m, q \nabla \cdot ] u\,dx -  \tilde \gamma\int_{\R^d} \pa^m u \cdot[\pa^m, q  \nabla] q\,dx.
\end{split}
\end{equation*}
We then estimate using Gagliardo--Nirenberg--Sobolev inequality that 
\begin{equation*}
\begin{split}
&\left|\frac{d}{dt}\left( \frac{1}{2}\int_{\R^d} |\rd^mu|^2 + |\rd^m q|^2\,dx \right)- \tilde{c}_K \int_{\R^d} \rd^m u\cdot \rd^m \nb \Lmb^{\alp-d}(q^N)  \,dx\right| \cr
&\quad \le C(\nrm{u}_{H^m} + \nrm{q}_{H^m}) \left( \int_{\R^d} |\rd^mu|^2 + |\rd^m q|^2 \,dx \right),
\end{split}
\end{equation*}
where we set $N = \frac{1}{\tilde{\gmm}}$. Let us write 
\begin{equation*}
\begin{split}
\int_{\R^d} \rd^m u\cdot \rd^m \nb \Lmb^{\alp-d}(q^{N}) \,dx &= N\int_{\R^d} \rd^m u\cdot \rd^m \Lmb^{\alp-d}(q^{N-1} \nb q)\,dx \cr
&= N\int_{\R^d} \rd^m u\cdot \Lmb^{\alp-d}(q^{N-1} \nb\rd^m  q) \,dx+R,
\end{split}
\end{equation*} where \begin{equation*}
\begin{split}
|R|\le C\sum_{0< |k|\le m} \left| \int_{\R^d} \rd^m u\cdot \Lmb^{\alp-d}( \rd^k(q^{N-1}) \nb\rd^{m-k}  q)\,dx\right| \le C \nrm{u}_{H^m} \nrm{q^{N-1}}_{H^{m}} \nrm{q}_{H^{m}}. 
\end{split}
\end{equation*} 
Thus we obtain
\begin{equation}\label{eq:concl-1}
\begin{split}
&\left| \frac{d}{dt}\left( \frac{1}{2}\int_{\R^d} |\rd^mu|^2 + |\rd^m q|^2 \,dx\right)- N\tilde{c}_K \int_{\R^d} \rd^m u\cdot \Lmb^{\alp-d}(q^{N-1} \nb\rd^m  q) \,dx\right|\\
 &\qquad \le C(\nrm{u}_{H^m} + \nrm{q}_{H^m}) (\nrm{ u}_{H^m}^2 + \nrm{q}_{H^m}(\nrm{q}_{H^m} + \nrm{q^{N-1}}_{H^m} )  ) .
\end{split}
\end{equation} 
On the other hand, we find
\begin{equation*}
\begin{split}
&\int_{\R^d} \rd^m u\cdot \Lmb^{\alp-d}(q^{N-1} \nb\rd^m  q) \,dx \cr
&\quad = \int_{\R^d} \rd^m u \cdot q^{N-1}\Lmb^{\alp-d}\nb\rd^m  q \,dx + \int_{\R^d} \pa^m u \cdot [\Lambda^{\alpha-d}, q^{N-1} \nabla] \pa^m q\,dx\cr
&\quad = \int_{\R^d}  (\rd^{m-1}\nb\cdot u)\, q^{N-1} \rd( \Lmb^{\alp-d} (\rd^mq)) \,dx + \int_{\R^d}  (\rd^{m-1}\nb\cdot u)\, \rd(q^{N-1}) ( \Lmb^{\alp-d} (\rd^mq)) \,dx\cr
&\qquad -  \int_{\R^d} \rd^m u \cdot \nabla(q^{N-1}) \Lmb^{\alp-d}\rd^m  q \,dx+ \int_{\R^d} \pa^m u \cdot [\Lambda^{\alpha-d}, q^{N-1} \nabla] \pa^m q\,dx.
\end{split}
\end{equation*}
Note that the second, third, and fourth terms on the right hand side of the above equality can be bounded from above by
\[
C\|u\|_{H^m} \|q^{N-1}\|_{H^m}\|q\|_{H^m}
\]
for some $C>0$. Thus we combine this with \eqref{eq:concl-1} to have 
\begin{equation}\label{eq:concl-1.5}
\begin{split}
&\left| \frac{d}{dt}\left( \frac{1}{2}\int_{\R^d} |\rd^mu|^2 + |\rd^m q|^2 \,dx\right)- N\tilde{c}_K \int_{\R^d}  (\rd^{m-1}\nb\cdot u)\, q^{N-1} \rd( \Lmb^{\alp-d} (\rd^mq)) \,dx\right|\\
 &\qquad \le C(\nrm{u}_{H^m} + \nrm{q}_{H^m}) (\nrm{ u}_{H^m}^2 + \nrm{q}_{H^m}(\nrm{q}_{H^m} + \nrm{q^{N-1}}_{H^m} )  ) .
\end{split}
\end{equation}

Let us now compute that 
\begin{equation*}
\begin{split}
&\frac{d}{dt} \frac{1}{2} \int_{\R^d} q^{N-2} |\Lmb^{\frac{\alp-d}{2}} \rd^mq|^2\,dx \cr
&\quad = \frac{1}{2} \int_{\R^d} \rd_t (q^{N-2}) |\Lmb^{\frac{\alp-d}{2}} \rd^mq|^2 \,dx+ \int_{\R^d} q^{N-2} \Lmb^{\frac{\alp-d}{2}} (\rd^mq) \Lmb^{\frac{\alp-d}{2}} \rd^m( -u\cdot\nb q - \tilde{\gmm} q\nb\cdot u )\,dx \\
&\quad =: I + J, 
\end{split}
\end{equation*} and $I$ is easily estimated as follows:
\begin{equation*}
\begin{split}
|I| \le C\left( \nrm{u}_{L^\infty} \nrm{\nb(q^{N-2})}_{L^\infty} + \nrm{\nb u}_{L^\infty} \nrm{q^{N-2}}_{L^\infty} \right) \int_{\R^d}  |\Lmb^{\frac{\alp-d}{2}} \rd^mq|^2 \,dx.
\end{split}
\end{equation*} Defining \begin{equation*}
\begin{split}
J_1 := -\tilde{\gmm}\int_{\R^d}  q^{N-2} \Lmb^{\frac{\alp-d}{2}} (\rd^mq) \Lmb^{\frac{\alp-d}{2}} (\rd (  q \, \rd^{m-1}\nb\cdot u ))\,dx,
\end{split}
\end{equation*} it is not difficult to see that \begin{equation*}
\begin{split}
\left| J - J_1 \right| \le C \nrm{q^{N-2}}_{L^\infty}   \nrm{u}_{H^m}\nrm{q}_{H^m}^2. 
\end{split}
\end{equation*} Integrating by parts, we obtain that \begin{equation*}
\begin{split}
J_1 = \tilde{\gmm} \int_{\R^d}  \Lmb^{\frac{\alp-d}{2}}\rd(q^{N-2} \Lmb^{\frac{\alp-d}{2}} (\rd^mq))   (  q \, \rd^{m-1}\nb\cdot u ) \,dx
\end{split}
\end{equation*} and writing $ \Lmb^{\frac{\alp-d}{2}}\rd(q^{N-2} \cdot ) = q^{N-2}  \Lmb^{\frac{\alp-d}{2}}\rd( \cdot ) + [ \Lmb^{\frac{\alp-d}{2}}\rd, q^{N-2} ](\cdot)$, \begin{equation*}
\begin{split}
J_1 &= \tilde{\gmm} \int_{\R^d} q^{N-1} \rd( \Lmb^{\alp-d} (\rd^mq))   \, \rd^{m-1}\nb\cdot u \,dx +  \tilde{\gmm} \int_{\R^d}  [ \Lmb^{\frac{\alp-d}{2}}\rd, q^{N-2} ](\Lmb^{\frac{\alp-d}{2}} (\rd^mq) )   \, (  q \, \rd^{m-1}\nb\cdot u )\,dx \cr
&=:J_{11}+J_{12}.
\end{split}
\end{equation*} Using Proposition \ref{prop:comm} with $0<\eps < \frac{d-\alp}{2} \wedge (m - \frac d2 - 1)$, we get 
\begin{equation*}
\begin{split}
\left|J_{12}\right| \le C \nrm{ q^{N-2} }_{H^m} \nrm{q}_{H^m} \nrm{q}_{L^\infty}\nrm{u}_{H^m}. 
\end{split}
\end{equation*} We obtain that \begin{equation}\label{eq:concl-2}
\begin{split}
&\left|\frac{d}{dt} \frac{1}{2} \int_{\R^d} q^{N-2} |\Lmb^{\frac{\alp-d}{2}} \rd^mq|^2 \,dx- \tilde{\gmm} \int_{\R^d} q^{N-1} \rd( \Lmb^{\alp-d} (\rd^mq))   \, \rd^{m-1}\nb\cdot u \,dx  \right| \cr
&\qquad \le C\nrm{ q^{N-2} }_{H^m} \nrm{u}_{H^m} \nrm{q}_{H^m}^2. 
\end{split}
\end{equation} Combining \eqref{eq:concl-1.5} and \eqref{eq:concl-2} implies
\begin{equation*}
\begin{split}
&\left| \frac{d}{dt}\left( \frac{1}{2}\int_{\R^d} |\rd^mu|^2 + |\rd^m q|^2 - N^2\tilde{c}_K q^{N-2} |\Lmb^{\frac{\alp-d}{2}} \rd^mq|^2 \,dx \right) \right|\\
&\qquad \le C(\nrm{u}_{H^m} + \nrm{q}_{H^m}) (1 + \nrm{ q^{N-2} }_{H^m}) (\nrm{ u}_{H^m}^2 + \nrm{q}_{H^m}(\nrm{q}_{H^m} + \nrm{q^{N-1}}_{H^m} )  ).
\end{split}
\end{equation*} Applying Lemma \ref{lem:power-Sobolev}, we find
\begin{equation*}
\begin{split}
\nrm{ q^{N-2} }^2_{H^{m}} & \le C\nrm{q^{N-3}}^2_{L^\infty} \nrm{q}_{H^m}^2 +C(1 + \nrm{\nb (\ln q)}_{L^\infty})^{2m} \nrm{q^{N-2}}_{L^2}^2 \\
& \le C \nrm{q}_{H^m}^{2(N-2)} + C(1 + \nrm{\nb (\ln q)}_{L^\infty})^{2m} \nrm{q}_{L^{2(N-2)}}^{2(N-2)} 
\end{split}
\end{equation*}  and then using the algebra property of $H^{m}$ with $\nrm{q}_{L^{2(N-2)}} \le C \nrm{q}_{H^m}$,\begin{equation*}
\begin{split}
\nrm{q^{N-1}}_{H^m} \le C \nrm{q^{N-2}}_{H^m} \nrm{q}_{H^m} \le C (1 + \nrm{\nb (\ln q)}_{L^\infty})^{ m} \nrm{q}_{H^m}^{ N-1}.
\end{split}
\end{equation*} This gives \begin{equation*}
\begin{split}
&\left| \frac{d}{dt}\left( \frac{1}{2}\int_{\R^d} |\rd^mu|^2 + |\rd^m q|^2 - N^2\tilde{c}_K q^{N-2} |\Lmb^{\frac{\alp-d}{2}} \rd^mq|^2 \,dx \right) \right|\\
&\qquad \le C (1 + \nrm{\nb (\ln q)}_{L^\infty})^{3m}  (1 + \nrm{ u}_{H^m}^2 + \nrm{q}_{H^m}^2 + \nrm{q}_{H^m}^{2(N-1)})^4.
\end{split}
\end{equation*} In the case $\tilde{c}_K>0$, we now estimate using the Plancherel theorem, with some $\xi_0>0$ to be determined, \begin{equation*}
\begin{split}
\int_{\R^d} q^{N-2} |\Lmb^{\frac{\alp-d}{2}} \rd^mq|^2 \,dx \le C \nrm{q}_{L^\infty}^{N-2}\left( \xi_0^{2m+\alp-d} \nrm{q}_{L^2}^2 + \xi_0^{\alp-d} \nrm{\rd^m q}_{L^2}^2   \right).
\end{split}
\end{equation*} Hence, by taking $\xi_0 \gg 1$ sufficiently large (note that it does not depend on $m$), we can guarantee that \begin{equation*}
\begin{split}
N^2\tilde{c}_K\int_{\R^d} q^{N-2} |\Lmb^{\frac{\alp-d}{2}} \rd^mq|^2  \,dx\le \frac{1}{10} \nrm{\rd^m q}_{L^2}^2  + C \nrm{q}_{L^\infty}^{M} \nrm{q}_{L^2}^2.
\end{split}
\end{equation*} Here $M = M(N,\alp,d)>0$ is some large power. In conclusion, we have the control \begin{equation*}
\begin{split}
\nrm{(u,q)}_{\dot{X}^m}^2 &:=  \int_{\R^d} |\rd^mu|^2 + |\rd^m q|^2 - N^2\tilde{c}_K q^{N-2} |\Lmb^{\frac{\alp-d}{2}} \rd^mq|^2\,dx + C \nrm{q}_{L^\infty}^{M} \nrm{q}_{L^2}^2\cr
& \ge \frac{1}{2} \int_{\R^d} |\rd^mu|^2 + |\rd^m q|^2 \,dx.
\end{split}
\end{equation*} This control is available \textit{a fortiori} for derivatives of order lower than $m$. Defining \begin{equation*}
\begin{split}
\nrm{(u,q)}_{ {X}^m}^2 := \nrm{u}_{L^2}^2+\nrm{q}_{L^2}^2 + \sum_{\ell=1}^{m} \nrm{(u,q)}_{\dot{X}^{\ell}}^2, 
\end{split}
\end{equation*} we have that \begin{equation*}
\begin{split}
\nrm{(u,q)}_{ {X}^m}^2 \ge \frac{1}{2}( \nrm{u}_{H^m}^2 + \nrm{q}_{H^m}^2 ).
\end{split}
\end{equation*} Hence, there exists some power $\tilde{M}>0$ such that \begin{equation*}
\begin{split}
\frac{d}{dt} \left( \nrm{(u,q)}_{X^m} + \nrm{\nb(\ln q)}_{L^\infty} \right) \lesssim \left(1+ \nrm{(u,q)}_{X^m} + \nrm{\nb(\ln q)}_{L^\infty} \right)^{\tilde{M}}.
\end{split}
\end{equation*} This finishes the proof of an a priori estimate. 

We now treat the case $\gmm =1$; the equation is now given by \eqref{eq:Euler-Riesz-pressure3}. Similarly as in the above proof, we  find that 
\begin{equation*}
\begin{split}
&\left| \frac{d}{dt}\left( \frac{1}{2}\int_{\R^d} |\rd^mu|^2 + |\rd^m q|^2 \,dx\right)- c_K \int_{\R^d} \rd^m u\cdot \pa^m \nabla \Lambda^{\alpha-d} (e^q) \,dx\right|\\
 &\qquad \le C(\nrm{u}_{H^m} + \nrm{q}_{H^m}) (\nrm{ u}_{H^m}^2 + \nrm{q}_{H^m}^2 ).
\end{split}
\end{equation*} 
Again by using a similar argument as in the proof of Theorem \ref{thm:lwp-pressure}, we obtain
\begin{equation*}
\begin{split}
&\left| \frac{d}{dt}\left( \frac{1}{2}\int_{\R^d} |\rd^mu|^2 + |\rd^m q|^2 \,dx\right)- c_K \int_{\R^d} e^q \rd( \Lmb^{\alp-d} (\rd^mq))   \, \rd^{m-1}\nb\cdot u \,dx\right|\\
 &\qquad \le C(\nrm{u}_{H^m} + \nrm{q}_{H^m}) (\nrm{ u}_{H^m}^2 + \nrm{q}_{H^m}^2 ) + Ce^{\nrm{q}_{H^m}} \nrm{u}_{H^m} \nrm{q}_{H^m}^2 .
\end{split}
\end{equation*} 
Similarly, we also estimate
\[
\left|\frac{d}{dt} \frac{1}{2} \int_{\R^d} e^q |\Lmb^{\frac{\alp-d}{2}} \rd^mq|^2 \,dx-   \int_{\R^d} e^q \rd( \Lmb^{\alp-d} (\rd^mq))   \, \rd^{m-1}\nb\cdot u \,dx  \right| \le Ce^{\nrm{q}_{H^m}} \nrm{u}_{H^m} \nrm{q}_{H^m}^2. 
\]
Combining the above two estimates implies
\begin{equation*}
\begin{split}
&\left| \frac{d}{dt}\left( \frac{1}{2}\int_{\R^d} |\rd^mu|^2 + |\rd^m q|^2 - c_Ke^q |\Lmb^{\frac{\alp-d}{2}} \rd^mq|^2 \,dx\right) \right|\\
 &\qquad \le C(\nrm{u}_{H^m} + \nrm{q}_{H^m}) (\nrm{ u}_{H^m}^2 + \nrm{q}_{H^m}^2 ) + Ce^{\nrm{q}_{H^m}} \nrm{u}_{H^m} \nrm{q}_{H^m}^2 .
\end{split}
\end{equation*}
On the other hand, we get
\[
|c_K| \int_{\R^d} e^q |\Lmb^{\frac{\alp-d}{2}} \rd^mq|^2 \,dx \leq |c_K|e^{\|q\|_{H^m}} \int_{\R^d}   |\Lmb^{\frac{\alp-d}{2}} \rd^mq|^2 \,dx \leq \frac{1}{10} \|\pa^m q\|_{L^2}^2 + Ce^{M\|q\|_{H^m}}\|q\|_{L^2}^2.
\]
This yields
\[
\nrm{(u,q)}_{\dot{Y}^m}^2 :=  \int_{\R^d} |\rd^mu|^2 + |\rd^m q|^2 - c_K e^q |\Lmb^{\frac{\alp-d}{2}} \rd^mq|^2\,dx + Ce^{M\|q\|_{H^m}} \nrm{q}_{L^2}^2  \ge \frac{1}{2} \int_{\R^d} |\rd^mu|^2 + |\rd^m q|^2 \,dx,
\]
and thus
\[
\nrm{(u,q)}_{Y^m}^2 := \nrm{u}_{L^2}^2+\nrm{q}_{L^2}^2 + \sum_{\ell=1}^{m} \nrm{(u,q)}_{\dot{Y}^{\ell}}^2 \geq \frac{1}{2}( \nrm{u}_{H^m}^2 + \nrm{q}_{H^m}^2 ).
\]
Hence we have 
\[
\frac{d}{dt}   \nrm{(u,q)}_{Y^m}  \lesssim e^{\tilde M\nrm{(u,q)}_{Y^m}}
\]
for some $\tilde M > 0$. This completes the proof of an a priori estimate.

Let us just briefly comment on the proof of uniqueness and existence. Uniqueness can be proved easily along the lines of the uniqueness proof for Theorem \ref{thm:pressureless-lwp}; essentially, one can perform a simple $L^2$ estimate for the difference of two hypothetical solutions associated with the same initial data. Existence of a solution follows from viscous approximations. \end{proof}

\begin{remark}\label{rem:pressure} Observe that in the statement of Theorem \ref{thm:lwp-pressure}, the range of $\gmm$ is restricted to $[1,\frac{5}{3}]$. In any of the following situations, such a restriction can be dropped:
	\begin{itemize}
		\item Either there exists $r>0$ such that $\rho_0(x)\ge r$ (in particular, when the domain is given by $\bbT^d$), or
		\item the potential is repulsive ($c_K \le 0$), or
		\item $-2 < \alp - d \le -1$.
	\end{itemize} In other words, the most difficult case is when $\rho_0$ decays at infinity, $c_K > 0$, and $-1 < \alp-d < 0$. Indeed, when $-2 < \alp - d \le -1$, the problematic term can be handled as follows: \begin{equation*}
\begin{split}
\left| \rd^m u \cdot \rd^m \nb \Lmb^{\alp-d} (q^N) \right| \le C \nrm{\rd^m u}_{L^2} \nrm{q^{N}}_{H^m},
\end{split}
\end{equation*} and it is easy to propagate $H^m$ regularity of $q^{N}$ with $u\in H^m$. Next, when $\rho_0 > 0$ admits a uniform lower bound, one can perform weighted energy estimates in the spirit of Theorem \ref{thm:pressureless-lwp}. To be more precise, one can close an a priori estimate in terms of \begin{equation*}
\begin{split}
\sum_{0\le\ell\le m}\left( \nrm{\rho^{-\frac{1}{2}} \rd^\ell \rho}_{L^2}^2  - c_K \nrm{\Lmb^{\frac{\alp-d}{2}} (\rd^\ell\rho)}_{L^2}^2 + \nrm{\rho^{\frac{1}{2}}\rd^\ell u}_{L^2}^2\right) 
\end{split}
\end{equation*} when $m$ is large enough. Here, $\nrm{\rho^{\frac{1}{2}}\rd^\ell u}_{L^2}$ is equivalent with $\nrm{\rd^\ell u}_{L^2}$ thanks to the lower bound of $\rho$. When $\rho$ decays at infinity, one can still try to propagate the modified quantity \begin{equation*}
\begin{split}
\sum_{0\le\ell\le m}\left( \nrm{\rho^{-1} \rd^\ell \rho}_{L^2}^2  - c_K \nrm{\rho^{-\frac{1}{2}}\Lmb^{\frac{\alp-d}{2}} (\rd^\ell\rho)}_{L^2}^2 + \nrm{ \rd^\ell u}_{L^2}^2\right) .
\end{split}
\end{equation*} This can be done with some additional work (there is a problem with the quantity $\nrm{\rho^{-1} \rd^\ell \rho}_{L^2}$ for small $\ell\ge 0$) when $c_K \le 0$. When $c_K>0$, it seems like that there is some serious problem in controlling the second quantity in terms of the weighted norm $\nrm{\rho^{-1} \rd^\ell \rho}_{L^2}$. 
\end{remark}

%
%
%
%
%

\section{Finite-time singularity formation}\label{sec:sing}

In this section, we analyze a finite-time singularity formation for the system \eqref{eq:Euler-second} in the presence of pressure.  For this, inspired by \cite{LPZ,S85,W14,Xi}, we introduce some physical quantities. $I = I(t)$ and $W=W(t)$ denote the moment of inertia and weighted momentum
\[
I:= \frac12\int_{\R^d} \rho|x|^2\,dx \quad \mbox{and} \quad W:= \int_{\R^d} \rho u \cdot x\,dx, 
\]
respectively. We also set a total energy $E_\gamma = E_\gamma(t)$
\[
E_\gamma := E_u + c_pE_{\rho,\gamma} -c_K E_K,
\]
where 
\[
E_u := \frac12\int_{\R^d} \rho |u|^2\,dx,  \quad E_K:=\frac{1}{2}\int_{\R^d} \rho \Lmb^{\alp-d}\rho\,dx,
\]
and
\[
E_{\rho,\gamma}  := \left\{ \begin{array}{ll}
\displaystyle \frac{1}{\gamma-1}\int_{\R^d} \rho^\gamma\,dx & \textrm{if $\gamma > 1$}\\[4mm]
\displaystyle \int_{\R^d} \rho \ln \rho\,dx& \textrm{if $\gamma=1$}
\end{array} \right..
\]

Before we proceed, we first define a solution space $\mathcal{Z}$ as follows.

\begin{definition}\label{def_sol} For a given $T>0$, we call $(\rho,u) \in \mathcal{Z}(T)$ if $(\rho,u)$ is a classical solution to the Cauchy problem \eqref{eq:Euler-second} on the time interval $[0,T]$ satisfying the following conditions of decay at far fields:
\[
\rho |u||x|^2 + (|u| + |x|) \lt( \rho |u|^2 + p(\rho) + \rho |\nabla \Lambda^{\alpha-d} \rho| \rt) \to 0 
\]
as $|x| \to +\infty$ for all $t \in [0,T]$.
\end{definition}
The decay conditions for solutions allow us to do the integration by parts in our estimates. We emphasize that all the estimates in this section are rigorous due to the existence theory in the previous section. More precisely, we can find a solution to the system \eqref{eq:Euler-second} which satisfies the assumptions in Theorem \ref{thm:lwp-pressure} and the decay conditions in Definition \eqref{def_sol}; note that our local well-posedness result covers $\rho$ which decays exponentially fast in space. Our strategy does not require that either the initial density $\rho_0$ contains vacuum or has compact support in any finite region. Moreover, it can be directly applied to the case of Coulombian interaction potentials, i.e. $\alpha = d-2$. 

We begin with the estimates of conservation laws for the system \eqref{eq:Euler-second}. Since its proofs are by now standard, we omit it here.

\begin{lemma}\label{lem_en} Let  $(\rho,u)$ be a solution to the system \eqref{eq:Euler-second} satisfying $(\rho,u) \in \mathcal{Z}(T)$.  Then we have
\[
\frac{d}{dt}\int_{\R^d} \rho\,dx = 0, \quad \frac{d}{dt}\int_{\R^d} \rho u\,dx = 0, \quad \mbox{and} \quad \frac{d}{dt} E_\gamma =0.
\]
\end{lemma}
In the lemma below, we provide some useful estimates giving relationships among the physical quantities defined as above.

\begin{lemma}\label{lem_aux} Let  $(\rho,u)$ be a solution to the system \eqref{eq:Euler-second} satisfying $(\rho,u) \in \mathcal{Z}(T)$. Then we have
\[
\frac{d}{dt}I(t) = W(t), \quad W(t)^2 \leq 4E_u(t) I(t),
\]
and
\[
\frac{d}{dt} W(t) = 2E_u(t) + c_p d(\gamma-1) E_{\rho,\gamma}(t) + c_p d\delta_{\gamma,1}\int_{\R^d} \rho(t,x)\,dx - \alpha c_KE_K(t),
\]
where $\delta_{\gamma,1}$ denotes the Kronecker delta. 
\end{lemma}
\begin{proof} We only prove the third assertion since the others are clear. A straightforward computation gives
\[
\frac{d}{dt}W = \int_{\R^d} \rho|u|^2\,dx + c_p d\int_{\R^d} \rho^\gamma\,dx + c_K \int_{\R^d} \rho x \cdot \nabla \Lmb^{\alp-d}\rho\,dx.
\]
On the other hand, the third term on the right hand side of the above can be estimated as 
$$\begin{aligned}
\int_{\R^d} \rho x \cdot \nabla \Lmb^{\alp-d}\rho\,dx &= -\alpha \iint_{\R^d \times \R^d} \rho(x) x \cdot \frac{x-y}{|x-y|^{\alpha+2}} \rho(y)\,dxdy\cr
&= \alpha \iint_{\R^d \times \R^d} \rho(y) y \cdot \frac{x-y}{|x-y|^{\alpha+2}} \rho(x)\,dxdy\cr
&=-\frac\alpha2 \iint_{\R^d \times \R^d} \rho(x) \frac{1}{|x-y|^\alpha} \rho(y)\,dxdy\cr
&=-\frac\alpha2\int_{\R^d} \rho \Lmb^{\alp-d}\rho\,dx\cr
&=-\alpha E_K.
\end{aligned}$$
This concludes the desired result.
\end{proof}

In the following two subsections, we present the finite-time singularity formation for our main system \eqref{eq:Euler-second} in the presence of pressure. For simplicity of presentation, without loss of generality, we set $c_p = 1$.

%
%
%
%
%

\subsection{Isentropic pressure case}
\subsubsection{Attractive case} In this part, we consider the attractive interaction case, i.e., $c_K >0$.  Without loss of generality, we set $c_K = 1$.

\begin{theorem}[Attractive case]\label{thm_att} Let $T>0$, $\gamma > 1$, and $(\rho,u)$ be a solution to the system \eqref{eq:Euler-second} satisfying $(\rho,u) \in \mathcal{Z}(T)$. Suppose that $\alpha > 0$ is large enough such that 
\bq\label{con_att1}
\alpha \geq \max\{2, d(\gamma-1)\}
\eq
and
\[
I(0), \ W(0), \ E_\gamma(0) < \infty.
\] 
If the initial data $(\rho_0, u_0)$ satisfies 
\bq\label{con_att2}
\frac12\int_{\R^d} \rho_0 |u_0|^2\,dx + \frac{1}{\gamma-1}\int_{\R^d} \rho_0^\gamma\,dx < \frac12\int_{\R^d} \rho_0 \Lmb^{\alp-d}\rho_0\,dx,
\eq
then the life-span $T$ of the solution is finite.
\end{theorem}
\begin{proof}[Proof of Theorem \ref{thm_att}] It follows from Lemma \ref{lem_aux} and \eqref{con_att1} that
$$\begin{aligned}
\frac{d}{dt} W(t) &= 2E_u(t) + d(\gamma-1) E_{\rho,\gamma}(t) - \alpha E_K(t)\cr
&\leq \max\{2, d(\gamma-1)\} \lt(E_u(t) + E_{\rho,\gamma}(t) \rt) - \alpha E_K(t)\cr
&= \max\{2, d(\gamma-1)\} \lt(E_\gamma(t) + E_K(t) \rt) - \alpha E_K(t)\cr
&=\max\{2, d(\gamma-1)\} E_\gamma(t) + \lt(\max\{2, d(\gamma-1)\} - \alpha\rt) E_K(t)\cr
&\leq\max\{2, d(\gamma-1)\} E_\gamma(t).
\end{aligned}$$
We then use Lemma \eqref{lem_en} to get
\[
\frac{d}{dt} W(t) \leq c_{d,\gamma} E_\gamma(t),
\]
where $c_{d,\gamma} := \max\{2, d(\gamma-1)\}$. We now integrate the above differential inequality over $[0,t]$ twice to find
\[
I(t) \leq I(0) + W(0) t + \frac12 c_{d,\gamma} E_\gamma(0) t^2.
\]
Since our assumption \eqref{con_att2} implies $E_0 < 0$, the right hand side is negative for sufficiently large $t$, while the left hand side remains nonnegative. This yields that the life-span $T$  of the solution should be finite. 
\end{proof}

\subsubsection{Repulsive case} We next consider the case $c_K < 0$. Again, for simplicity, we set $c_K = -1$. Note that in this case we cannot use the strategy used in the previous subsection since the initial total energy is always nonnegative. In order to overcome this, we need to get the more detailed information on the upper bound of $I$. More precisely, we obtain from Lemma \ref{lem_aux} that 
\[
\frac{d}{dt} W(t) = 2E_u(t) + d(\gamma-1) E_{\rho,\gamma}(t) + \alpha E_K(t) \leq \max\{2, d(\gamma-1), \alpha \} E_\gamma(t) = \max\{2, d(\gamma-1), \alpha \} E_\gamma(0).
\]
This implies
\bq\label{upp_I}
I(t) \leq I(0) + W(0) t + \frac{c_{d,\gamma,\alpha}}{2} E_\gamma(0) t^2,
\eq
where $c_{d,\gamma,\alpha} :=  \max\{2, d(\gamma-1), \alpha \} > 0$. 

Using the time-growth estimate of $I$, we show the lower bound estimate of $E_{\rho,\gamma}$ in the lemma below.
\begin{lemma}\label{lem_I}Let $(\rho,u)$ be a solution to the system \eqref{eq:Euler-second} satisfying $(\rho,u) \in \mathcal{Z}(T)$. Then we have
\[
E_{\rho,\gamma}(t) \geq \frac{c_0}{\lt(I(0) + W(0) t + \frac{c_{d,\gamma,\alpha}}{2} E_\gamma(0) t^2\rt)^{\frac{d(\gamma-1)}{2}}},
\]
where $c_0 > 0$ is given by
\[
c_0 := \lt(\frac{\pi^{d/2}}{\Gamma(d/2+1)} \rt)^{1-\gamma} \frac{\|\rho_0\|_{L^1}^{\frac{(d+2)\gamma-d}{2}}}{2^{\frac{(d+2)\gamma-d}{2}(\gamma-1)}}.
\]
Here $\Gamma$ is the gamma function. 
\end{lemma}
\begin{proof}
For any $R>0$, we estimate
\[
\int_{\R^d} \rho\,dx = \lt(\int_{|x| \leq R} + \int_{|x| \geq R} \rt) \rho\,dx \leq |B(0,R)|^{1-\frac1\gamma}\lt(\int_{\R^d} \rho^\gamma\,dx \rt)^{\frac1\gamma} + \frac{1}{R^2}\int_{\R^d} \rho|x|^2\,dx,
\]
where $B(0,R) := \{x \in \R^d : |x| \leq R\}$ and $|A|$ denotes the Lebesgue measure of a set $A$ in $\R^d$. We then choose $R$ so that the right hand side of the above inequality is minimized, i.e.,
\[
R = \lt(\frac{\int_{\R^d} \rho|x|^2\,dx}{\|\rho\|_{L^\gamma} |B(0,1)|^{1-\frac1\gamma}} \rt)^{\frac{\gamma}{(d+2)\gamma - d}}.
\] 
This yields
\[
E_{\rho,\gamma}(t) \geq \frac{c_0}{I(t)^{\frac{d(\gamma-1)}{2}}}.
\]
We finally combine this with \eqref{upp_I} to complete the proof.
\end{proof}
Our next goal is to estimate the upper bound on $E_{\rho,\gamma}$. For this, we set
\bq\label{eq_J}
J(t) = E_\gamma(t)(t+1)^2 - W(t) (t+1) + I(t).
\eq
Then by the second assertion in Lemma \ref{lem_aux} we get
\bq\label{ineq_J}
J(t) \geq (t+1)^2 (E_{\rho,\gamma}(t) + E_K(t)).
\eq
On the other hand, differentiating \eqref{eq_J} with respect to $t$ yields
\[
J'(t) = -(t+1) W'(t) + 2(t+1)E_\gamma(t) = (t+1)(2 - d(\gamma-1)) E_{\rho,\gamma}(t) + (t+1)(2-\alpha) E_K(t)
\]
due to Lemmas \ref{lem_en} and \ref{lem_aux}. We now assume 
\[
2-\alpha \leq 2 - d(\gamma-1), \quad \mbox{i.e.,} \quad d(\gamma-1) \leq \alpha
\]
and
\[
1 < \gamma \leq 1 + \frac{2}{d}.
\]
Then we have
\[
J'(t) \leq (t+1)(2 - d(\gamma-1)) (E_{\rho,\gamma}(t) + E_K(t)),
\]
and this together with \eqref{ineq_J} gives
\[
J'(t) \leq \frac{2 - d(\gamma-1)}{(t+1)} J(t).
\]
We finally solve the above differential inequality to have the following lemma.
\begin{lemma}\label{lem_J}Let $(\rho,u)$ be a solution to the system \eqref{eq:Euler-second} satisfying $(\rho,u) \in \mathcal{Z}(T)$. Assume
\[
1 < \gamma \leq 1 + \frac{2}{d} \quad \mbox{and} \quad d(\gamma-1) \leq \alpha.
\]
Then we have
\[
J(t) \leq \frac{J(0)}{(t+1)^{d(\gamma-1) - 2}}.
\]
\end{lemma}
We now state our result on the finite-time singularity formation in the repulsive interaction potential case. 
\begin{theorem}[Repulsive case]\label{thm_rep}  Let $T>0$ and $(\rho,u)$ be a solution to the system \eqref{eq:Euler-second} satisfying $(\rho,u) \in \mathcal{Z}(T)$. Suppose that $\gamma$ and $\alpha > 0$ satisfy
\[
1 < \gamma \leq 1 + \frac{2}{d} \quad \mbox{and} \quad \alpha \geq d(\gamma-1),
\]
respectively. Moreover, we assume
\[
I(0), \ W(0), \ E_\gamma(0), \ \mbox{and} \ \|\rho_0\|_{L^1} < \infty.
\] 
If the initial data $(\rho_0, u_0)$ satisfies 
\bq\label{con_rep}
J(0) = E_\gamma(0) - W(0) + I(0)  <\frac{2c_0}{c_{d,\gamma,\alpha} E_\gamma(0)},
\eq
where $c_0$ is given in Lemma \ref{lem_I}, then the life-span $T$ of the solution is finite.
\end{theorem}
\begin{proof} We first estimate the lower bound on $J$. Since $E_K \geq 0$, we get from \eqref{ineq_J} 
\[
J(t) \geq (t+1)^2 E_{\rho,\gamma}(t).
\]
We further use Lemma \ref{lem_I} to obtain
\[
J(t) \geq \frac{c_0(t+1)^2 }{\lt( I(0) + W(0) t + \frac{c_{d,\gamma,\alpha}}{2} E_\gamma(0) t^2\rt)^{\frac{d(\gamma-1)}{2}}}.
\]
This and Lemma \ref{lem_J} yield
\[
\frac{J(0)}{(t+1)^{d(\gamma-1) - 2}} \geq \frac{c_0(t+1)^2 }{\lt( I(0) + W(0) t + \frac{c_{d,\gamma,\alpha}}{2} E_\gamma(0) t^2\rt)^{\frac{d(\gamma-1)}{2}}}.
\]
Note that as $t$ tends to infinity the above inequality implies
\[
J(0) \geq \frac{2c_0}{c_{d,\gamma,\alpha} E_\gamma(0)},
\]
and this contradicts \eqref{con_rep}. Hence the life-span $T$ of the solution should be finite.
\end{proof}

%
%
%
%
%

\subsection{Isothermal pressure case}
In this subsection, we deal with the isothermal pressure law, i.e., $p(\rho) = \rho$ in the system \eqref{eq:Euler-second}. 

As a direct application of Lemma \ref{lem_aux}, we first find
\bq\label{est_iso}
\frac{d^2}{dt^2} I(t) = \frac{d}{dt}W(t) = 2E_u + d\int_{\R^d} \rho\,dx - \alpha c_K E_K = 2E_1  + (2-\alpha)c_K E_K + d\|\rho\|_{L^1} - 2E_{\rho,1},
\eq
where $E_1$ is the total energy for the system \eqref{eq:Euler-second} with $\gamma=1$, i.e., $E_1= E_u + E_{\rho,1}  - c_K E_K$. It is worth noticing that the last term on right hand side of the above inequality does not have the definite sign. This requires the control of the negative part of $E_{\rho,1}$. For this, we provide the auxiliary lemma below.
\begin{lemma}\label{lem_iso} For a given $\epsilon > 0$, there exists a $C_0>0$ independent of $\rho$ such that 
\[
-\epsilon\int_{\R^d} \rho (\ln \rho) \chi_{0 \leq \rho \leq 1}\,dx \leq \frac12\int_{\R^d}  \rho |x|^2\,dx + C_\epsilon.
\]
\end{lemma}
\begin{proof}Note that the following holds for any $s, \sigma \geq 0$:
$$\begin{aligned}
-s\ln(s) \chi_{0 \leq s \leq 1} &= -s\ln(s) \chi_{e^{-\sigma} \leq s \leq 1} - s\ln(s) \chi_{e^{-\sigma} \geq s}\cr
&\leq s\sigma + C\sqrt s \chi_{e^{-\sigma} \geq s}\cr
&\leq s\sigma + Ce^{-\sigma^2/2}
\end{aligned}$$
for some $C>0$ independent of $s$ and $\sigma$. This together with taking $s = \rho$ and $\sigma = |x|^2/{2\epsilon}$ asserts
\[
-\epsilon\int_{\R^d} \rho (\ln \rho) \chi_{0 \leq \rho \leq 1}\,dx \leq \frac12\int_{\R^d}  \rho |x|^2\,dx + C_\epsilon
\]
for some $C_\epsilon > 0$.
\end{proof}

We now state our result on the finite-time singularity formation for the isothermal pressure case. 

\begin{theorem}\label{thm_iso} Let $T>0$, $\gamma = 1$, and $(\rho,u)$ be a solution to the system \eqref{eq:Euler-second} satisfying $(\rho,u) \in \mathcal{Z}(T)$. Suppose 
\[
I(0), \ W(0), \ E_\gamma(0), \ \mbox{and} \ \|\rho_0\|_{L^1} < \infty.
\] 
If $c_K = 1$, i.e., in the attractive case, then we assume that $\alpha \geq 2$ and the initial data $(\rho_0, u_0)$ satisfies 
\[
W(0) + I(0) + 2E_1(0) + d\|\rho_0\|_{L^1} + C_2 < 0.
\]
On the other hand, if $c_K = -1$, i.e., in the repulsive case, we suppose that the initial data $(\rho_0, u_0)$ satisfies 
\[
W(0) + I(0) +\max\{2,\alpha \} E_1(0) + d\|\rho_0\|_{L^1} + C_{\max\{2,\alpha \}} < 0
\]
for $\alpha \geq 0$ and 
\[
W(0) + I(0) + 2E_1(0) + d\|\rho_0\|_{L^1} + C_2 < 0
\]
for $\alpha < 0$. Here $C_2$ and $C_{\max\{2,\alpha \}}$ are positive constants given as in Lemma \ref{lem_iso}. Then the life-span $T$ of the solution is finite.
\end{theorem}
\begin{remark}Note that either $W_0$ or $E_0$ should be chosen to be negative to make our assumptions in Theorem \ref{thm_iso} valid. 
\end{remark}
\begin{proof}[Proof of Theorem \ref{thm_iso}] We first derive a second-order differential inequality for $I(t)$. If $c_K = 1$, i.e., in the attractive case, then by choosing $\alpha$ large enough so that $\alpha \geq 2$, we obtain from \eqref{est_iso} that
\[
\frac{d^2}{dt^2} I(t) \leq 2E_1(t) + d\|\rho_0\|_{L^1}   - 2E_{\rho,1}(t).
\]
Then this together with Lemma \ref{lem_iso} and total energy estimate gives
$$\begin{aligned}
\frac{d^2}{dt^2} I(t) &\leq 2E_1(0)   + d\|\rho_0\|_{L^1} - 2\int_{\R^d} \rho(t,x) (\ln \rho(t,x)) \chi_{0 \leq \rho(t,x) \leq 1}\,dx\cr
&\leq 2E_1(0) + d\|\rho_0\|_{L^1} + C_2 + I(t),
\end{aligned}$$
where $C_2$ is given in Lemma \ref{lem_iso}. 

For the repulsive case, i.e., $c_K = -1$, we estimate 
\[
\frac{d^2}{dt^2} I(t) = 2E_u + d\|\rho\|_{L^1} + \alpha E_K \leq \max\{2,\alpha \} E_1(t) + d\|\rho\|_{L^1} - \max\{2,\alpha \} E_{\rho,1}(t)
\]
for $\alpha \geq 0$. Then similarly as before we deduce 
\[
\frac{d^2}{dt^2} I(t)  \leq \max\{2,\alpha \} E_1(0) + d\|\rho_0\|_{L^1} + C_{\max\{2,\alpha \}} + I(t).
\]
On the other hand, if $\alpha < 0$, then we get
\[
\frac{d^2}{dt^2} I(t) \leq 2E_u + d\|\rho\|_{L^1} \leq 2E_1(t) + d\|\rho_0\|_{L^1}   - 2E_{\rho,1}(t) \leq 2E_1(0) + d\|\rho_0\|_{L^1} + C_2 + I(t).
\]
Thus for both cases we have the following form of second-order differential inequality:
\bq\label{est_I2}
\frac{d^2}{dt^2} I(t) \leq \tilde C + I(t), 
\eq
with the initial data $I(0)$ and $I'(0) = W(0)$, where $\tilde C>0$ is given by
\[
\tilde C = \left\{ \begin{array}{ll}
2E_1(0) + d\|\rho_0\|_{L^1} + C_2 & \textrm{if $c_K = 1$}\\[2mm]
\max\{2,\alpha \} E_1(0) + d\|\rho_0\|_{L^1} + C_{\max\{2,\alpha \}} & \textrm{if $c_K = -1$ and $\alpha \geq 0$}\\[2mm]
2E_1(0) + d\|\rho_0\|_{L^1} + C_2 & \textrm{if $c_K = -1$ and $\alpha < 0$}
\end{array} \right..
\]
In order to solve the second-order differential inequality \eqref{est_I2}, by introducing 
\[
\tilde I(t) := \frac{d}{dt}I(t) - I(t),
\]
we first reduce it to the first-order differential inequality:
\[
\frac{d}{dt} \tilde I(t) + \tilde I(t) \leq \tilde C,
\]
subject to the initial data $\tilde I(0) = W(0) - I(0)$. We then apply Gr\"onwall's lemma twice to obtain
\[
I(t) \leq (I(0) + \tilde C)e^{t} + \frac12 \lt(e^t - e^{-t} \rt)(W(0) - I(0) - \tilde C) - \tilde C.
\]
Note that our assumptions on the initial data give $W(0) + I(0) + \tilde C < 0$ and thus the right hand side is negative for $t>0$ large enough. This is contradiction that $I(t)$ is always nonnegative. Hence the life-span $T$ of solution should be finite. 
\end{proof}

%
%
%
%
%

\section{Conclusion}

In this paper, we initiated the analysis of the Euler--Riesz system. We established a local well-posedness theory for the Euler--Riesz system (with pressure)  and pressureless and repulsive Euler--Riesz system. This well-posedness result makes the mean-field limit \cite[Appendix A]{Sepre} from the Newton's second-order particle system \eqref{par_sys} to the pressureless Euler--Riesz system \eqref{eq:Euler-second} fully rigorous. Our framework for well-posedness allows the fluid density to decay fast at infinity, and it covers the Euler--Poison system. Finally, we investigated the problem of finite-time singularity formation for classical solutions for the Euler--Riesz system with $\gamma \geq 1$. We provided sufficient conditions on the initial data, $\gamma$, $\alpha$, and $d$ which lead to the finite-time breakdown of smoothness of solutions. Unfortunately, the strategy used for Theorem \ref{thm_rep} cannot be applied to the pressureless case since the internal energy, which appears due to the presence of pressure, plays a crucial role in the analysis of blow-up of solutions. We would like to stress that the finite-time singularity for the pressureless Euler--Poisson with repulsive forces is still a challenging open problem. We would like to investigate this question in the future.

\section*{Acknowledgement} We thank Bongsuk Kwon for helpful conversations regarding the Euler--Poisson system. 
YPC has been supported by NRF grant (No. 2017R1C1B2012918), POSCO Science Fellowship of POSCO TJ Park Foundation, and Yonsei University Research Fund of 2019-22-021. IJJ has been supported  by a KIAS Individual Grant MG066202 at Korea Institute for Advanced Study, the Science Fellowship of POSCO TJ Park Foundation, and the National Research Foundation of Korea grant (No. 2019R1F1A1058486).

\bibliographystyle{amsplain}
\bibliography{Euler}

\end{document}